\documentclass[11pt]{article}
\usepackage{hyperref}
\usepackage{amsfonts,amsmath,amssymb,mathtools,bbm,amsthm,color}
\usepackage{cleveref}
\usepackage{lineno}
\modulolinenumbers[10]
\usepackage{geometry}
\usepackage{tabstackengine}
\usepackage{amscd,amsbsy,epsf}
\usepackage{enumitem}
\usepackage{titlesec}
\usepackage{graphicx}
\usepackage[normalem]{ulem}
\usepackage{tikz,sansmath,booktabs}
\usetikzlibrary{arrows,backgrounds}
\usetikzlibrary{shapes.geometric, arrows}
\usetikzlibrary[calc]
\usetikzlibrary{positioning}
\usetikzlibrary{fit}
\usetikzlibrary{arrows.meta}
\tikzset{block_device_small/.style={draw, thick, text width=1.5cm, minimum height=1.0cm, align=left,fill=cyan},   
}

\tikzset{block_device/.style={draw, thick, text width=4.5cm, minimum height=1.5cm, fill={cyan},  align=center}}
\tikzset{block_laptop/.style={draw, thick, text width=4cm, minimum height=1.5cm, fill={pink},  align=center}}
\tikzstyle{arrow} = [thick,-Stealth]
\tikzstyle{line} = [thick,-]
\tikzstyle{square-node} = [rectangle, rounded corners, minimum width=3cm, minimum height=1cm,text centered, draw=black, , fill=red!30]

\tikzset{
    *|/.style={
        to path={
            (perpendicular cs: horizontal line through={(\tikztostart)},
                                 vertical line through={(\tikztotarget)})
            -- (\tikztotarget) \tikztonodes
        }
    }
}
\usepackage{algorithm,algorithmicx}
\usepackage{algpseudocode}
\usepackage{soul,wrapfig,float}

\def\calD {{\mathcal D}}

\newcommand{\mD}{{\mathcal D}}

\newcommand{\R}{\mathbb{R}}

\newcommand{\bmu}{\boldsymbol \mu}

\newcommand{\bc}{\mathbf c}

 \newcommand{\bx}{\mathbf x}

\newcommand{\cC}{\mathcal C}  
 \newcommand{\cF}{\mathcal F}
\newcommand{\cG}{\mathcal G}

\usepackage{multirow}

\begin{document}
\title{GPT-PINN: Generative Pre-Trained Physics-Informed Neural Networks toward non-intrusive Meta-learning of parametric PDEs\footnote{This work was partially supported by National Science Foundation grant DMS-2208277, by Air Force Office of Scientific Research grant FA9550-23-1-0037, and by the UMass Dartmouth Marine and UnderSea Technology (MUST) Research Program made possible via an Office of Naval Research grant N00014-20-1-2849. The code of GPT-PINN is available at \href{https://github.com/skoohy/GPT-PINN}{https://github.com/skoohy/GPT-PINN}. Authors' address: Department of Mathematics, University of Massachusetts Dartmouth, 285 Old Westport Rd, North Dartmouth, 02747, MA USA.
}}
\author{Yanlai Chen\footnote{Corresponding Author. Email address: {\tt yanlai.chen@umassd.edu}} , \, Shawn Koohy\footnote{Email address: {\tt skoohy@umassd.edu}}}

\date{}
\maketitle
\begin{abstract}
Physics-Informed Neural Network (PINN) has proven itself a powerful tool to obtain the numerical solutions of nonlinear partial differential equations (PDEs) leveraging the expressivity of deep neural networks and the computing power of modern heterogeneous hardware. However, its training is still time-consuming, especially in the multi-query and real-time simulation settings, and its parameterization often overly excessive. In this paper, we propose the Generative Pre-Trained PINN (GPT-PINN) to mitigate both challenges in the setting of parametric PDEs. GPT-PINN represents a brand-new meta-learning paradigm for parametric systems. As a network of networks, its outer-/meta-network is hyper-reduced with only one hidden layer having significantly reduced number of neurons. Moreover, its activation function at each hidden neuron is a  (full) PINN pre-trained at a judiciously selected system configuration. 
The meta-network adaptively ``learns'' the parametric dependence of the system and ``grows'' this hidden layer one neuron at a time. In the end, by encompassing a very small number of networks trained at this set of adaptively-selected parameter values, the meta-network is capable of generating surrogate solutions for the parametric system across the entire parameter domain accurately and efficiently. 
\end{abstract}
\section{Introduction}

The need to efficiently and accurately understand the behavior of the system under
variation of a {large number of underlying parameters} is ubiquitous in \emph{many query} type of applications e.g. uncertainty quantification, (Bayesian) inverse problems, data assimilation or optimal control/design. The  
parameters of interest may include material properties, wave
frequencies, uncertainties, boundary conditions, the shape of the domain, etc. 
A rigorous study of the behavior of the system and its dependence on the parameters  requires 
thousands, perhaps millions of simulations of the underlying partial differential equations (PDE).
Each accurate and robust simulation of the underlying complex physical phenomena is often time consuming,
and the massively repeated simulations needed become computationally challenging, if not entirely untenable, when using traditional numerical methods. 
Two techniques stand out in addressing this challenge, the more traditional and rigorous reduced order modeling and the more nascent deep neural networks.

The reduced basis method (RBM) \cite{patera_reduced_2007,quarteroni_reduced_2016,HesthavenRozzaStammBook, Rozza_Huynh_Patera,Maday_Patera_Rovas,Haasdonk2017Review}, a projection-based model order reduction approach \cite{BennerGugercinWillcox2015}, belongs to the first category. It was developed to generate a computational emulator for parameterized problems whose error compared to the full problem is certifiable; this rigorous accuracy guarantee is a relatively unique ability among reduced order model algorithms. Once generated, the RBM emulator, using the results of the original method at carefully preselected parameter values, can typically compute an accurate solution with orders-of-magnitude less computational cost than the original method. 
This is achieved through an offline-online decomposition, where the parameter values are selected and reduced solution space 
constructed in an offline  preparation/training phase via a greedy algorithm, allowing the rapid online computation of the surrogate solution for any parameter values.

Deep learning algorithms are increasingly popular in this context as well.  
Using data generated by many queries of the underlying system, one can train a deep neural network (DNN, a highly nonlinear function composed of layers of parameterized affine linear functions and simple nonlinear operations) to provide a surrogate for the parameter to solution map. Physics-informed neural networks (PINNs), popularized by \cite{raissi2019physics}, adopt DNNs to represent the approximate solutions of PDEs. 
Unlike typical data-driven deep learning methods that do not build in physical understanding of the problem,  PINNs incorporate a strong physics prior (i.e. PDEs) that constrains the output of the DNN. 
The key advantages of PINNs, over traditional numerical solvers, include that they are able to solve the PDE without discretizing the problem domain, that they define a function over the entire continuous spatial-temporal domain, and that they can rely on automatic differentiation \cite{baydin2017automatic,paszke2017automatic} toward residual minimization. 
Thanks also to the enormous advances in computational capabilities in recent years \cite{abadi2016tensorflow, revels2016forward}, 
PINNs have emerged as an increasingly popular alternative to traditional numerical methods for PDEs.

However, issues of PINNs remain \cite{wang2022respecting}. 
Among them, vanilla PINNs are usually significantly slower than the classic numerical methods due to the training of the, what is usually substantially parameterized, neural network. The main purpose of this paper is to use strategies inspired by the classical and mathematically rigorous RBM techniques to significantly shrink the size of PINNs and accelerate solving parametric PDEs with PINNs. 
Just like RBM, the proposed solvers have an initial investment cost. However, they are capable of providing significant computational savings in problems where a PDE must be solved repeatedly or in real-time thanks to the fact that their marginal cost is of orders of magnitude lower than that of each PINN solve. 

The jump from the vanilla PINN to the proposed Generative Pre-Trained PINN (GPT-PINN) parallels that from the traditional Finite Element Method (FEM) to RBM. To the best of our knowledge, it represents a first-of-its-kind meta-learning approach for parametric systems. Its infrastructure is a network of networks. The inner networks are the full PINNs. Its outer-/meta-network is hyper-reduced, in comparison to the inner networks, with only one hidden layer where the  inner networks are pre-trained and serve as activation functions. 
The meta-network adaptively ``learns'' the parametric dependence of the system and ``grows'' this hidden layer one neuron/network at a time. In the end, by encompassing a very small number of networks trained at this set of adaptively-selected parameter values, the meta-network is capable of generating surrogate solutions for the parametric system across the entire parameter domain accurately and efficiently, with a cost independent of the size of the full PINN. The design of network architecture represent the first main novelty of the paper. To the best of our knowledge, this is the first time whole (pre-trained) networks are used as the activation functions of another network. The adoption of the training loss of the meta-network as an error indicator, inspired by the residual-based error estimation for traditional numerical solvers such as FEM, represents the second main novelty.

The rest of the paper is organized as follows. In Section \ref{sec:bg}, we review the RBM and PINN. In Section \ref{sec:gpt_pinn}, we detail our design of GPT-PINN, and also remark on recent efforts about accelerating PINNs in the parametric PDE setting. We present numerical results on three parametric PDEs in Section \ref{sec:numerics} demonstrating the accuracy and efficiency of the proposed GPT-PINN. Finally, concluding remarks are given in Section \ref{sec:conclusion}.

\section{Background}
\label{sec:bg}
\subsection{Reduced Basis Method}
\label{sec:rbm}

RBM is a linear reduction method that has been a popular option for rigorously and efficiently simulating parametric PDEs. Its hallmark feature is a greedy algorithm embedded in an offline-online decomposition procedure. The offline (i.e. training) stage is devoted to a judicious and error estimate-driven exploration of the parameter-induced solution manifold. It  selects a number of representative parameter values via a mathematically rigorous greedy algorithm \cite{BinevCohenDahmenDevorePetrovaWojtaszczyk}. During the online stage,   a {\em reduced} solution is sought in the terminal surrogate space for each unseen parameter value. Moreover, unlike other reduction techniques (e.g. proper orthogonal decomposition (POD)-based approaches), the number of full order inquiries RBM takes offline is minimum, i.e. equal to the dimension of the surrogate space. To demonstrate the main ideas, we consider a
generic parameterized PDE as follows, 
\begin{equation}
\label{eq:general-pde}
  \mathcal{F}(u;\bx,\bmu) = f, \quad  x \in \Omega \subseteq \mathbb{R}^d.
  \end{equation}
Here $\mathcal{F}$ encodes a  differential operator parameterized via $\bmu \in \calD \subset \mathbb{R}^{d_s}$ together with necessary boundary and initial conditions. The parameter can be equation coefficients, initial values, source terms, or uncertainties in the PDE for the tasks of the uncertainty quantification, etc. $\mathcal{F}$ can depend on the solution and its (space- and time-) derivatives of various orders. 
We assume that we have available a numerical solution $u(\bx; \bmu) \in X_h$  obtained by a high fidelity solver, called Full Order Model (FOM) and denoted as FOM$({\bmu}, X_h)$, and $X_h$ is the discrete approximation space the numerical solution $u$ belongs to.

A large number of queries of $u(\cdot;\bmu)$ can be prohibitively expensive because the  FOM$({\bmu}, X_h)$ has to be called many times. Model order reduction (MOR) aims to mitigate this cost by building efficient surrogates. One idea is to study the map 
\[
\bmu \mapsto u(\cdot, \bmu) \in {X_h}
\]
and devise an algorithm to compute an approximation $u_N(\cdot, \mu)$ {from  an $N$-dimensional subspace $X_N$ of $X_h$,} 
such that 
\[
u_N(\cdot, \bmu) \approx u(\cdot, \bmu) 
\mbox{ for all } \bmu \in \calD
\]
This reduced order model (ROM) formulation at a given $\bmu$ is denoted by  ROM$({\bmu, X_N}),$ and is much cheaper to solve than  FOM$({\bmu}, X_h)$ and can be conducted during the Online stage.

\begin{algorithm}[H]
    \caption{Classical RBM for parametric PDE \eqref{eq:general-pde}: Offline stage}
    \label{alg:rbm_greedy}
    {\bf Input: }{A (random or given) $\bmu^1$, training set $\Xi \subset \calD$.}\\
    {\bf Initialization:} Solve FOM($\bmu^1, X_h$) and set $X_1 = \mbox{span}\left\{u(\cdot; \bmu_1)\right\}$, $n=2$.
\begin{algorithmic}[1]
\While{{\em stopping criteria not met,}}
\State Solve ROM($\bmu, X_{n-1}$) for all $\bmu \in \Xi$ and compute error indicators $\Delta_{n-1}(\bmu)$.
 \State Choose $\bmu^n = \displaystyle
  \mbox{\rm arg}\hspace*{-1pt}\max_{\bmu\in{\Xi}}
 {\Delta_{n-1}(\bmu)}$.
\State  Solve FOM($\bmu_n, X_h$) and update $X_n = X_{n-1} \bigoplus \{ u(\cdot; \bmu_n)\}$.
\State Set $n \leftarrow n+1$.
\EndWhile
\end{algorithmic} 
 {\bf{Output:}~}  Reduced basis set $X_N$, with {$N$ being the terminal index.}
\end{algorithm}
The success of RBM relies on the assumption that $u(\cdot; \calD)$ has small \textit{Kolmogorov $N$-width} \cite{pinkus_n-widths_1985}, defined as
\[
  d_N \left[ u\left(\cdot; \calD\right) \right] \coloneqq \inf_{\substack{X_N \subset {X_h} \\ \dim X_N = N}}\;\; \sup_{\bmu \in \calD}\;\; \inf_{v \in X_N} \left\| u(\cdot, \bmu) - v \right\|_X.
\]
A small $d_N$ means that the solution to \cref{eq:general-pde} for any $\bmu$ can be well-approximated from $X_N$ that represents the outer infimum above. The identification of a near-infimizing subspace $X_N$ is one of the central goals of RBM, and is obtained in the so-called Offline stage. RBM uses a greedy algorithm to find such $X_N.$ The main ingredients are presented in Algorithm \ref{alg:rbm_greedy}. The method explores the training parameter set {$\Xi\subset \mD$} guided by   an   error estimate {or an efficient and effective error indicator} $\Delta_n (\bmu)$  and intelligently choosing the parameter ensemble $\left\{\bmu^n\right\}_{n=1}^N$ so that 
\begin{equation}
X_N \coloneqq \mathrm{span}\left\{u(\cdot; \bmu^n)\right\}_{n=1}^N, \mbox{ and } u_N(\cdot, \bmu) = \sum_{n=1}^N c_n(\bmu) u(\cdot, \bmu^n). 
\label{eq:rbm_soln}
\end{equation}
An offline-online decomposed framework is key to realize the speedup.

Equipped with this robust and tested greedy algorithm, physics-informed reduced solver, rigorous error analysis, and certifiable convergence guarantees, 
RBM algorithms have become the go-to option for efficiently simulating parametric PDEs  and established in the modern scientific computing toolbox \cite{Prudhomme_Rovas_Veroy_Maday_Patera_Turinici,patera_reduced_2007,quarteroni_reduced_2016,HesthavenRozzaStammBook} and have benefited from voluminous research with theoretical and algorithmic refinement \cite{Noor_Peters,Peterson,Barrett_Reddien,Nagy,Rozza_Huynh_Patera,Maday_Patera_Rovas}. 
One particular such development was the empirical error indicator of the L1-based RBM by Chen and his collaborators \cite{JiangChenNarayan2019} where $\Delta_{n-1}(\mu)$ was taken to be $\left\lVert \bc(\mu) \right\rVert_1$.
Here $\bc(\mu)$ is the coefficient vector of $u_N(\cdot, \mu)$ under the basis $\left\{u(\cdot; \mu_n)\right\}_{n=1}^N$ and $\left\lVert \cdot \right\rVert_1$ represents the $\ell^1$-norm. As shown in \cite{JiangChenNarayan2019,ChenGottliebJiMaday2021}, $\bc(\mu)$ represents a Lagrange interpolation basis in the parameter space  implying that  the indicator $\Delta_n$ represents the corresponding Lebesgue constant. The L1 strategy to select the parameter samples then controls the growth of the Lebesgue constants and hence is key toward accurate interpolation. 
This strategy, ``free'' to compute albeit not as traditionally rigorous, inspires the greedy algorithm of our GPT-PINN, to be detailed in Section \ref{sec:gpt_pinn}.

\subsection{Deep neural networks}

Deep neural networks (DNN) have seen tremendous success recently when serving as universal approximators to the solution function (or certain quantity of interest (QoI) / observable)   \cite{lagaris_artificial_1998,perdikaris_nonlinear_2017,e_deep_2018,khoo_solving_2018,raissi_deep_2018, Cy1,YAR1, MSDR}. First proposed in \cite{lagaris_artificial_1998} on an underlying collocation approach, it has been successfully used recently in different contexts. See \cite{Kar1,Kar2, Kar3, JR1,e_deep_2018,E2,E3} and references therein. 
For a nonparametrized version (e.g. \cref{eq:general-pde} with a fixed parameter value), we search for a neural network $\Psi_{\mathsf{NN}}(\bx)$ which maps the coordinate $ \bx \in \mathbb{R}^d $ to a surrogate of the solution, that is $\Psi_{\mathsf{NN}} (\bx) \approx u (\bx) $.

Specifically, for an input vector $\bx$, a feedforward neural network maps it to an output, via layers of ``neurons''  with layer $k$ corresponding to an affine-linear map $C_k$ composed with scalar non-linear activation functions $\sigma$ \cite{DLbook}. That is,
\[
\Psi_{\mathsf{NN}}^\theta(\bx) = C_K \circ\sigma \circ C_{K-1}\ldots \ldots \circ \sigma \circ C_1(\bx).
\]
A justifiably popular choice is the \emph{ReLU} activation $\sigma(z) = \max(z,0)$ that is understood as component-wise operation when $z$ is a vector.  For any $1 \leq k \leq K$, we define
\[
C_k z_k = W_k z_k + b_k, \quad {\rm for} ~ W_k \in \R^{d_{k+1} \times d_k}, z_k \in \R^{d_k}, b_k \in \R^{d_{k+1}}.
\]
To be consistent with the input-output dimension, we set $d_1 = d$ and $d_K = 1$. 
We concatenate the tunable weights and biases for our network and denote them as 
\[
\theta \coloneqq \{W_k, b_k\},
\quad \forall~ 1 \leq k \leq K.
\]
We have $\theta \in \Theta \subset \R^M$ with $M \coloneqq \sum\limits_{k=1}^{K-1} (d_k +1) d_{k+1}$. We denote this network by
\begin{equation}
\label{eq:nn_notation}
\mathsf{NN}(d_1,d_2, \cdots, d_K).
\end{equation}
Learning $\Psi_{\mathsf{NN}}^\theta(\bx)$ then amounts to generating training data and determining the weights and biases $\theta$ by optimizing a loss function using this data.

\subsection{Physics-Informed Neural Network}

\label{sec:PINN}

We define our problem on the spatial domain $\Omega \subset \mathbb{R}^{d}$ with boundary $\partial \Omega$, and consider time-dependent PDEs with order of time-derivative $k = 1$ or $2$.
\begin{align}
\label{eq:pinn-pde}
\begin{split}
\frac{\partial^k}{\partial t^k}u(\bx, t) + \mathcal{F}\left[ u(\bx, t)\right] &= 0 \textnormal{  \hspace{1.67cm}   } \bx \in \Omega,\textnormal{  \hspace{1.0cm}   } t \in [0, T], \\
\mathcal{G}(u)(\bx, t) &= 0 \textnormal{  \hspace{1.67cm}   } \bx \in \partial\Omega,\textnormal{  \hspace{0.76cm}   } t \in [0, T], \\
u(\bx, 0) &= {u}_0(\bx) \textnormal{  \hspace{1.0cm}   } \bx \in \Omega.
\end{split}
\end{align}
Here $\mathcal{F}$ is a differential operator as defined in Section \ref{sec:rbm} and $\mathcal{G}$ denotes a boundary operator. The goal of a PINN is to identify an approximate solution $u(\bx,t)$ via a neural network $\Psi^\theta_{\mathsf{NN}}(\bx,t)$. 
Learning $\theta \in \mathbb{R}^M$ requires defining a loss function whose minimum $\theta^*$ leads to  $\Psi^{\theta^*}_{\mathsf{NN}}$ approximating the solution to the PDE over the problem domain. PINN defines this loss as a sum of three parts, an integral of the local residual of the differential equation over the problem domain, that over the boundary, and the deviation from the given initial condition,
\begin{align*}
  \mathcal{J}(u) = & \int_{\Omega} \left\lVert \frac{\partial^k}{\partial t^k}u(\bx, t) + \mathcal{F}(u)(\bx, t)\right\rVert^2_2 + \left\lVert u(\bx, 0) - {u}_0(\bx) \right\rVert _2^2\, dx +
   \int_{\partial\Omega} \left\lVert \mathcal{G}(u)(\bx, t)\right\rVert _2^2 \, dx.
\end{align*}
During training, we sample collocation points in certain fashion from the PDE space domain $\Omega$, space-time domain~$\Omega \times (0, T)$, and boundary $\partial \Omega \times [0,T]$, $\cC_o \subset \Omega \times (0,T)$ and~${\cC_\partial \subset \partial \Omega \times [0, T]}$ and ${\cC_i \subset \Omega}$, and use them to form an approximation of the true loss. 
\begin{align}
\label{eq:pinn_training_loss_time}
\begin{split}
    \mathcal{L}_{\text{PINN}}(\Psi^\theta_{\mathsf{NN}}) = & \frac{1}{|\cC_o|} \sum_{(\bx,t) \in \cC_o}\left\lVert \frac{\partial^k}{\partial t^k}(\Psi^\theta_{\mathsf{NN}})(\bx, t) + \mathcal{F} (\Psi^\theta_{\mathsf{NN}})(\bx, t) \right\rVert_2^2  +\\ 
     & \frac{1}{|\cC_\partial|}  \sum_{{(\bx,t) \in \mathcal{\cC_\partial}} } \left\lVert\mathcal{G}(\Psi^\theta_{\mathsf{NN}})(\bx, t)\right\rVert_2^2 + \frac{1}{|\cC_i|} \sum_{\bx \in \cC_i} \left\lVert \Psi^\theta_{\mathsf{NN}}(\bx, 0) - {u}_0(\bx) \right\rVert _2^2.
\end{split}
\end{align}
When the training converges, we expect that  $\mathcal{L}_{\text{PINN}}(\Psi^\theta_{\mathsf{NN}})$ should be nearly zero.

\section{The GPT-PINN framework}
\label{sec:gpt_pinn}

\begin{figure}[htbp]
    \centering
    \includegraphics[width=0.99\textwidth]{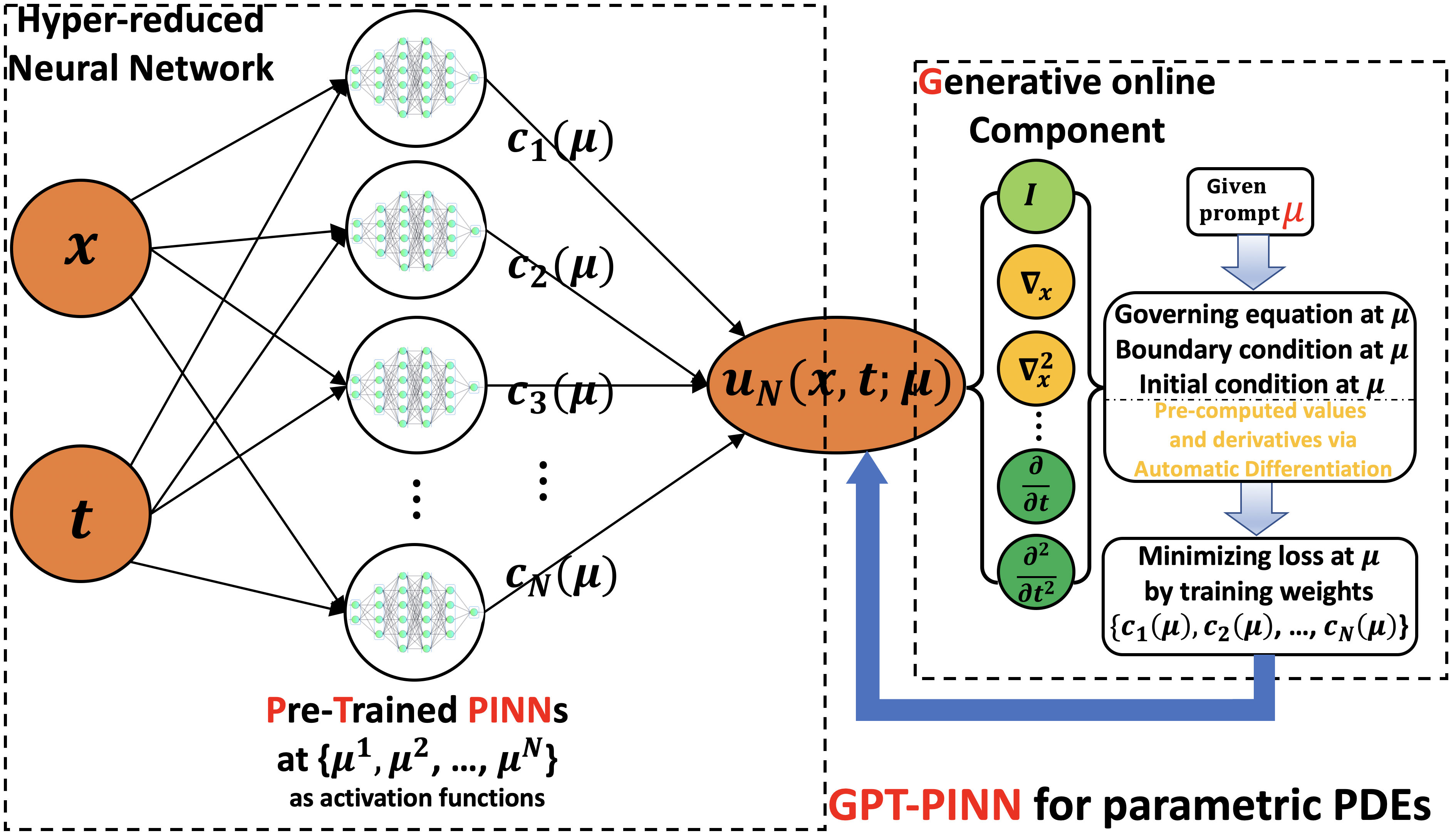}
    \caption{The GPT-PINN architecture. A hyper-reduced network adaptively embedding pre-trained PINNs at the nodes of its sole hidden layer. It then allows a quick online generation of a surrogate solution at any given parameter value.}
    \label{fig:gptpinn_diagram}
\end{figure}

Inspired by the RBM formulation \cref{eq:rbm_soln}, we design the GPT-PINN. Its two components and design philosophy are depicted in Figure \ref{fig:gptpinn_diagram}. As a hyper-reduced feedforward neural network $\mathsf{NN}(2,n,1)$ with $1 \le n \le N$ (see \cref{eq:nn_notation} for the notation), we denoted it by $\mathsf{NN}^{\rm r}(2,n,1)$.  A key feature is that it has customized  activation function in the neurons of its sole hidden layer. These activation functions are nothing but the pre-trained PINNs for the corresponding PDEs instantiated by the parameter values $\{\bmu^1, \bmu^2, \cdots, \bmu^n\}$ chosen by a greedy algorithm that is specifically tailored for PINNs but inspired by the classical one adopted by RBM in Algorithm \ref{alg:rbm_greedy}. The design of network architecture represents the first main novelty of the paper. To the best of our knowledge, this is the first time a whole (pre-trained) network is used as the activation function of one neuron.

\subsection{The online solver of GPT-PINN}

\label{sec:onlinesolver}

We first present the online solver, i.e. the training of the reduced network $\mathsf{NN}^{\rm r}(2, n, 1)$, for any given $\bmu$. With the next subsection detailing how we ``grow'' the GPT-PINN offline from $\mathsf{NN}^{\rm r}(2, n, 1)$ to $\mathsf{NN}^{\rm r}(2, n+1, 1)$, we have a strategy of adaptively generating the terminal GPT-PINN, $\mathsf{NN}^{\rm r}(2, N, 1)$. 
Indeed, given the simplicity of the reduced network, to train the weights $\{c_1(\bmu), \cdots, c_n(\bmu)\}$, no backpropagation is needed. The reason is that the loss function, similar to \cref{eq:pinn_training_loss_time}, is a simple function containing directly and explicitly $\{c_1(\bmu), \cdots, c_n(\bmu)\}$ thanks to the reduced network structure of GPT-PINN. 
In fact, we denote by $\Psi_{\mathsf{NN}}^{\theta^i}(x,t)$ the PINN approximation of the PDE solution when $\bmu = \bmu^i$. Given that $u_n(\bx,t; \bmu) \approx \sum_{i=1}^n c_i(\bmu) \Psi_{\mathsf{NN}}^{\theta^i}(x,t)$, we can calculate the GPT-PINN loss as a function of the weights $\bc(\bmu)$ as follows.

\begin{align}
\label{eq:loss-online}
\begin{split}
\mathcal{L}_{\text{PINN}}^{\text{GPT}}&(\bc(\bmu)) =  \frac{1}{|\cC_o^r|} \sum_{(\bx,t) \in \cC_o}\left\lVert \frac{\partial^k}{\partial t^k}\left (\sum_{i=1}^n c_i(\bmu)\Psi^{\theta^i}_{\mathsf{NN}}\right)(\bx, t) + \mathcal{F} \left(\sum_{i=1}^n c_i(\bmu)\Psi^{\theta^i}_{\mathsf{NN}}\right)(\bx, t) \right\rVert_2^2  +\\ 
     & \frac{1}{|\cC_\partial^r|}  \sum_{{(\bx,t) \in \mathcal{\cC_\partial}} } \left\lVert\mathcal{G}\left(\sum_{i=1}^n c_i(\bmu)\Psi^{\theta^i}_{\mathsf{NN}}\right)(\bx, t)\right\rVert_2^2 + \frac{1}{|\cC_i^r|} \sum_{\bx \in \cC_i} \left\lVert \sum_{i=1}^n c_i(\bmu) \Psi^{\theta^i}_{\mathsf{NN}}(\bx, 0) - {u}_0(\bx) \right\rVert _2^2.
\end{split}
\end{align}

The online collocation sets 
$\cC_o^r \subset \Omega \times (0,T)$, ${\cC_\partial^r \subset \partial \Omega \times [0, T]}$ and ${\cC_i^r \subset \Omega}$ are used, similar to \cref{eq:pinn_training_loss_time}, to generate an approximation of the true loss. They are taken to be the same as their full PINN counterparts $\cC_o, \cC_\partial, \cC_i$ in this paper but we note that they can be fully independent. The training of $\mathsf{NN}^{\rm r}(2, n, 1)$ is then simply
\begin{equation}
\label{eq:c-update}
    \bc \leftarrow \bc - \delta_r \nabla_\bc \mathcal{L}_{\text{PINN}}^{\text{GPT}}(\bc)
\end{equation}
Here $\bc = \left(c_1(\bmu), \cdots, c_n(\bmu) \right)^T$ and $\delta_r$ is the online learning rate. The detailed calculations of \cref{eq:loss-online} and \cref{eq:c-update} are given in  \ref{sec:appendix} for the first numerical example. Those for the other examples are very similar and thus omitted. We make the following three remarks to conclude the online solver.

\noindent {\bf 1. Precomputation for fast training of $\mathsf{NN}^{\rm r}(2, n, 1)$:} Due to the linearity of the derivative operations and the collocation nature of loss function, a significant amount of calculations of \cref{eq:loss-online} can be precomputed and stored. These include the function values and all (spatial and time) derivatives involved in the operators $\cF$ and $\cG$ of the PDE \cref{eq:pinn-pde}:
\begin{equation}
\label{eq:precompute-quantities}
    \Psi^{\theta^i}_{\mathsf{NN}}(\cC), \, \frac{\partial^k}{\partial t^k}\left (\Psi^{\theta^i}_{\mathsf{NN}}\right)(\cC)\, (k=1 \mbox{ or } 2), \, \nabla_\bx^\ell \Psi^{\theta^i}_{\mathsf{NN}}(\cC) \, (\ell = 1, 2, \cdots) \mbox{ for } \cC = \cC_o^r, \cC_\partial^r, \cC_i^r.
\end{equation}
Once these are precomputed, updating $\bc$ according to \cref{eq:c-update} is very efficient. It can even be made independent of $|\cC|$.

\noindent {\bf 2. Non-intrusiveness of GPT-PINN:} It is clear that, once the quantities of \cref{eq:precompute-quantities} are extracted from the full PINN, the online training of $\mathsf{NN}^{\rm r}(2, n, 1)$ is independent of the full PINN. GPT-PINN is therefore non-intrusive of the Full Order Model. One manifestation of this property is that, as shown in our third numerical example, the full PINN can be adaptive while the reduced PINN may not be.

\noindent {\bf 3. The error indication of $\mathsf{NN}^{\rm r}(2, n, 1)$.} One prominent feature of RBM is its {\em a posteriori} error estimators/indicators which guides the generation of the reduced solution space and certifies the accuracy of the surrogate solution. Inspired by this classical design, we introduce the following quantity that measures how accurate $\mathsf{NN}^{\rm r}(2, n, 1)$ is in generating a surrogate network at a new parameter $\bmu$.
\begin{equation}
\label{eq:ee_gpt_pinn}
    \Delta_{\mathsf{NN}}^r(\bc(\bmu)) \triangleq \mathcal{L}_{\text{PINN}}^{\text{GPT}}(\bc(\bmu)).
\end{equation}
We remark that this quantity is essentially free since it is readily available when we train $\mathsf{NN}^{\rm r}(2, n, 1)$ according to \cref{eq:c-update}. The adoption of the training loss of the meta-network as an error indicator, inspired by the residual-based error estimation for traditional numerical solvers such as FEM, represents the second main novelty of this paper.

\subsection{Training the reduced network GPT-PINN: the greedy algorithm}

\begin{algorithm}[htbp]
    \caption{GPT\_PINN for parametric PDE: Offline stage}
    \label{alg:gptpinn_greedy}
    {\bf Input: }{A (random or given) $\bmu^1$, training set $\Xi_{\rm train} \subset \calD$, full PINN.}
\begin{algorithmic}[1]
\State Train a full PINN at $\bmu^1$ to obtain $\Psi_{\mathsf{NN}}^{\theta^1}$. Precompute quantities necessary for $\nabla_\bc \mathcal{L}_{\text{PINN}}^{\text{GPT}}$ at collocation nodes $\cC_o^r$, $\cC_\partial^r$, and $\cC_i^r$, see \cref{eq:precompute-quantities}. Set $n=2$.
\While{{\em stopping criteria not met,}}
\State Train $\mathsf{NN}^{\rm r}(2, n-1, 1)$ at $\bmu$ for all $\bmu \in \Xi_{\rm train}$ and record the indicator $\Delta_{\mathsf{NN}}^r(\bmu)$.
 \State Choose $\bmu^n = \displaystyle
  \mbox{\rm arg}\hspace*{-1pt}\max_{\bmu\in{\Xi_{\rm train}}}
 {\Delta_{\mathsf{NN}}^r(\bmu)}$.
\State  Train a full PINN at $\bmu^n$ to obtain $\Psi_{\mathsf{NN}}^{\theta^n}$. Precompute quantities necessary for $\nabla_\bc \mathcal{L}_{\text{PINN}}^{\text{GPT}}$ at collocation nodes $\cC_o^r$, $\cC_\partial^r$, and $\cC_i^r$, see \cref{eq:precompute-quantities}.
\State Update the GPT\_PINN by adding a neuron to the hidden layer to construct $\mathsf{NN}^{\rm r}(2, n, 1)$.
\State Set $n \leftarrow n+1$.
\EndWhile
\end{algorithmic} 
 {\bf{Output:}~}  GPT\_PINN $\mathsf{NN}^{\rm r}(2, N, 1)$, with {$N$ being the terminal index.}
\end{algorithm}

\begin{figure}[htbp]
    \centering
    \includegraphics[width=0.9\textwidth]{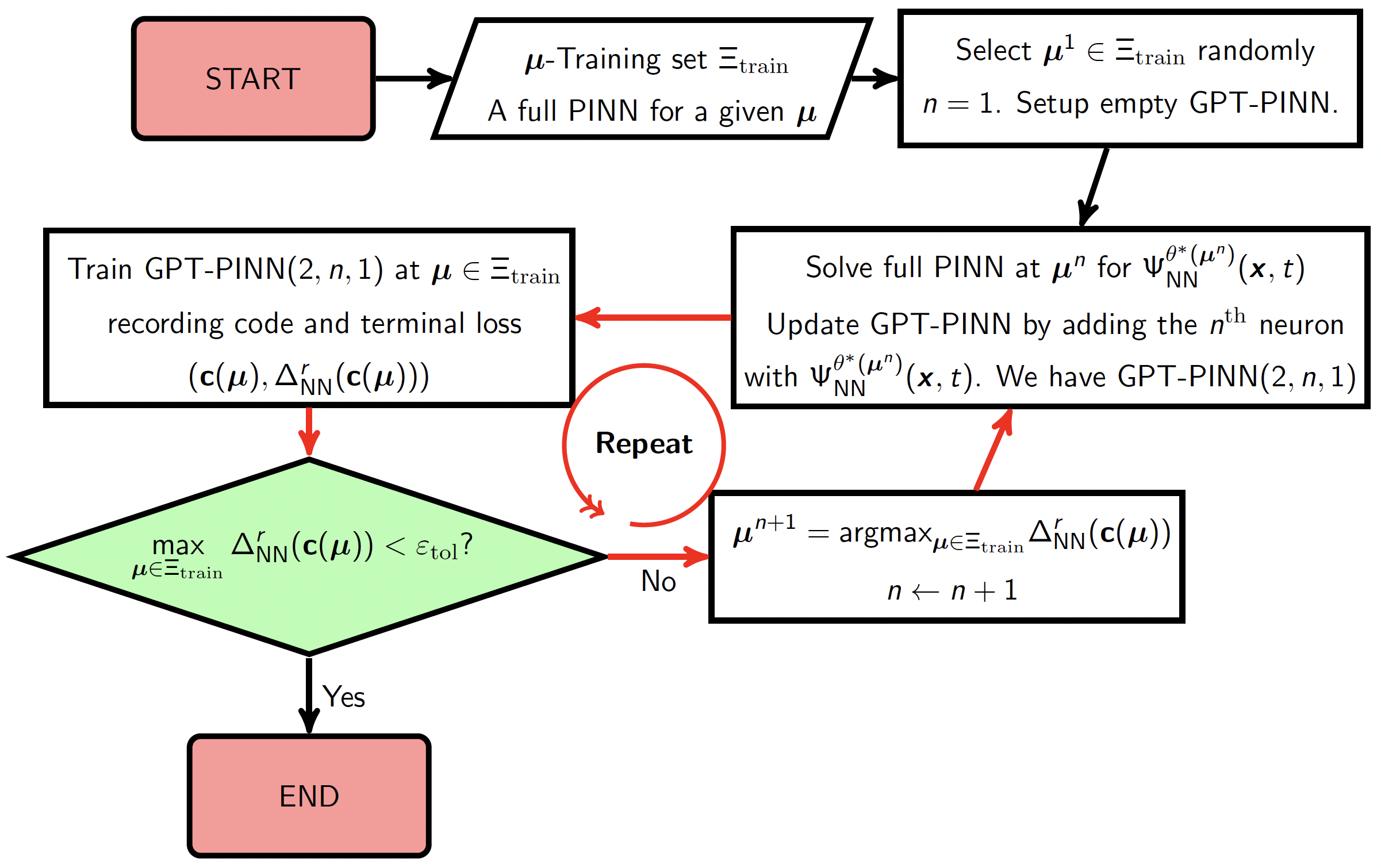}
    \caption{Flowchart of the GPT-PINN Offline training stage.}
    \label{fig:flowchart}
\end{figure}
With the online solver described in Section \ref{sec:onlinesolver}, we are ready to present our greedy algorithm. Its main steps are outlined in Algorithm \ref{alg:gptpinn_greedy} with its flowchart provided in Figure \ref{fig:flowchart}.  
The meta-network adaptively ``learns'' the parametric dependence of the system and ``grows'' its sole hidden layer one neuron/network at a time in the following fashion. We first randomly select, in the discretized parameter domain $\Xi_{\rm train}$, one parameter value $\bmu^1$ and train the associated (highly accurate) PINN $\Psi_{\mathsf{NN}}^{\theta^1}$. The algorithm then decides how to ``grow'' its meta-network by scanning the entire discrete parameter space $\Xi_{\rm train}$ and, for each parameter value, training this reduced network (of 1 hidden layer with 1 neuron $\Psi_{\mathsf{NN}}^{\theta^1}$). As it scans, it records an error indicator $\Delta_{\mathsf{NN}}^r(\bc(\bmu))$. The next parameter value $\bmu^2$ is the one generating the largest error indicator. The algorithm then proceeds by training a full PINN at $\bmu^2$ and therefore grows its hidden layer into two neurons with customized (but pre-trained) activation functions $\Psi_{\mathsf{NN}}^{\theta^1}$ and $\Psi_{\mathsf{NN}}^{\theta^2}$. This process is repeated until the stopping criteria is met which can be either that the error indicator is sufficiently small or a pre-selected size of the reduced network is met. At every step, we select the parameter value that is approximated most badly by the current meta-network. We end by presenting how we initialize the weights $\bc(\bmu)$ when we train $\mathsf{NN}^{\rm r}(2, n-1, 1)$ on Line 3 of Algorithm \ref{alg:gptpinn_greedy}. They are initialized by a linear interpolation of up to $2^{d_s}$ closest neighbors of $\bmu$ within the chosen parameter values $\{\bmu^1, \dots, \bmu^N\}$. Recall that $d_s$ is the dimension of the parameter domain.

\subsection{Related work}

The last two to three years have witnessed an increasing level of interest toward metalearning of (parameterized or unparameterized) PDEs due to the need of repeated simulations and the remarkable success of PINNs in its original form or adaptive ones. Here we mention a few representative ones and point out how our method differentiates from theirs.

\noindent {\bf Metalearning via PINN parameters.} In \cite{penwarden2021physics}, the authors adopt statistical (e.g. regression) and numerical (e.g. RBF/spline interpolation) methods to build a surrogate for the map from the PDE parameter $\bmu$ to the PINN parameter (weights and biases, $\theta$). They are shown to be superior than MAML \cite{finn2017model} for parameterized PDEs which was shown to outperform LEAP \cite{flennerhag2018transferring} in \cite{Qin_MetaPDE2022}. Both are general-purpose meta-learning methods. However, the online solver (i.e. regression or interpolation) of \cite{penwarden2021physics} ignores the physics (i.e. PDE). The method assumes that the $\bmu$-variation of the PINN weights and biases is analogous to that of the PDE solution.

\noindent {\bf DeepONet.} Aiming to learn nonlinear operators, a DeepONet \cite{lu2021learning} consists of two sub-networks, a branch net for encoding the input function (e.g source/control term, as opposed to PDE coefficients) at a fixed number of sensors, and a trunk net for encoding the locations for the output functions. It does not build in the physics represented by the dynamical system or PDE for a new input. Moreover, it is relatively data-intense by having to scan the entire input function space such as Gaussian random field or orthogonal polynomial space.

\noindent {\bf Metalearning loss functions.} Authors of \cite{PSAROS2022111121} concern the definition of the PINN loss functions. While it is in the parameterized PDE setting, the focus is a gradient-based approach to discover, during the offline stage, better PINN loss functions which are parameterized by e.g. the weights of each term in the composite objective function. The end goal is therefore improved PINN performance e.g. at unseen PDE parameters, due to the learned loss function configuration.

\noindent {\bf Metalearning initialization.} In \cite{Zhong_2023}, the authors study the use of a meta network, across the parameter domain of a 1-D arc model of plasma simulations, to better initialize the PINN at a new task (i.e. parameter value).

\noindent {\bf MetaNO.} The recent meta-learning approach for transferring knowledge
between neural operators \cite{zhang2023metano} aims to transfer the learned network parameters $\theta(\bmu)$ between different $\bmu$ with only the first layer being retrained. Its resulting surrogate is fully data-driven, i.e. with no physics built in for a new value $\bmu$.

\noindent {\bf PRNN.} The physics-reinforced neural network approach \cite{ChenWangHesthavenZhang2021} builds the map $\bmu \mapsto \bc(\bmu)$ via regression (i.e. no physics during the online evaluation for a new $\bmu$) although PDE residuals were considered during the supervised learning of the map via labelled data.

\begin{figure}[htbp]
    \centering
    \includegraphics[width=0.5\textwidth]{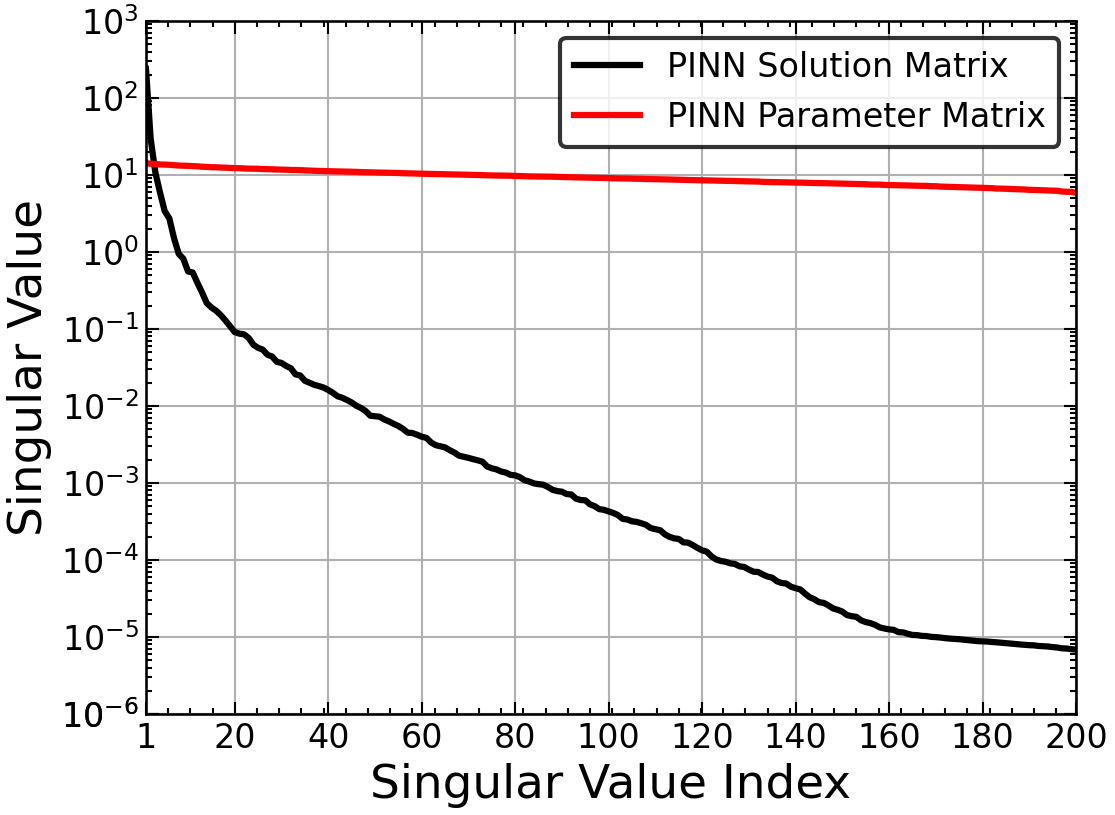}
    \caption{A motivating example showing that the solution matrix of a parametric PDE $\{u(\cdot, \bmu^n)\}_{n=1}^{200}$ exhibits fast decay in its singular values (indicating fast decay of the  Kolmogorov N-width of the solution manifold) while the network weights and biases manifold $\{\theta(\bmu^n)\}_{n=1}^{200}$ does not.}
    \label{fig:svd_example}
\end{figure}

Our proposed GPT-PINN exploits the $\bmu$-variation of the PDE solution directly which may feature a Kolmogorov N-width friendlier to MOR approaches, see Figure \ref{fig:svd_example}, than the weights and biases. This is, in part, because that the weights and biases lie in a (much) higher dimensional space. Moreover, the meta-network of our approach, being a PINN itself, has physics automatically built in in the same fashion as the underlying PINNs. Lastly, our approach provides a surrogate solution to the unseen parameter values in addition to a better initialization transferred from the sampled PINNs. 
Most importantly, our proposed GPT-PINN embodies prior knowledge that is mathematically rigorous and PDE-pertinent into the network architecture. This produces strong inductive bias that usually leads to good generalization.

\section{Numerical results}
\label{sec:numerics}

In this section, we present numerical results of the GPT-PINN applied to three families of equations, the Klein-Gordon equation, the Burgers' equation, and the Allen-Cahn equation. All simulations are run on a desktop with AMD Ryzen 7 2700X CPU clocked at 4.0GHz, an NVIDIA GeForce RTX 2060 SUPER GPU, and 32 GB of memory. Python version 3.9.12 was used along with common numerical packages and machine learning frameworks such as NumPy (v1.23.4), PyTorch (v1.11.0), TensorFlow (v2.10.0), and for GPU support CUDA v11.6 was installed. Previous literature \cite{wight2020solving,MATTEY2022114474,mcclenny2020self,xu2022numerical} has shown common difficulties in the use of baseline (non-adaptive) PINNs for the approximation of the Allen-Cahn equations. 
We have therefore adopted the Self-Adaptive PINNs (SA-PINNs) formulated by \cite{mcclenny2020self} in \cref{SEC-AC} to acquire accurate approximations by the full PINN, later used by the GPT-PINN. The tuned hyperparameters of the full PINNs include the activation functions and the learning rates. 
The code for all these examples are published on GitHub at \href{https://github.com/skoohy/GPT-PINN}{https://github.com/skoohy/GPT-PINN}.
Throughout the experiments of \cref{SEC-KG,SEC-B,SEC-AC}, we calculate and report various losses and errors. They are defined in Table \ref{tab:reportedquantities}.
\begin{table}[htbp]
  \begin{center}
  {
    \renewcommand{\tabcolsep}{0.4cm}
    {
    \renewcommand{\arraystretch}{1.5}
    \renewcommand{\tabcolsep}{12pt}
    \begin{tabular}{lr}
      \toprule
      Largest Loss &  The worst-case training loss ${\displaystyle \max_{\bmu \in \Xi_{\rm train}} \mathcal{L}_{\text{PINN}}^{\text{GPT}}(\bc(\bmu))}$ of $\mathsf{NN}^{\rm r}(2, n, 1)$ 
\\
      Terminal Losses & Training losses $\mathcal{L}_{\text{PINN}}^{\text{GPT}}(\bc(\bmu))$ of $\mathsf{NN}^{\rm r}(2, n, 1)$ (for more statistics)\\
      Largest Error & The worst-case testing error ${\displaystyle \max_{\bmu \in \Xi_{\rm test}} \frac{\left\lVert\mathsf{NN}^{\rm r}(2, n, 1)(x,t) - \Psi^{\theta(\bmu)}_{\mathsf{NN}}(x,t) \right\rVert_2}{\left\lVert\Psi^{\theta(\bmu)}_{\mathsf{NN}}(x,t)\right\rVert_2}}$\\
      Terminal Errors & Testing errors ${\displaystyle \frac{\left\lVert\mathsf{NN}^{\rm r}(2, n, 1)(x,t) - \Psi^{\theta(\bmu)}_{\mathsf{NN}}(x,t) \right\rVert_2}{\left\lVert\Psi^{\theta(\bmu)}_{\mathsf{NN}}(x,t)\right\rVert_2}}$ (for more statistics)\\
      Point-wise Error & Absolute error for a given $\bmu$, $\left |\mathsf{NN}^{\rm r}(2, N, 1)(x,t)- \Psi^{\theta}_{\mathsf{NN}}(x,t)\right|$\\
        \bottomrule
    \end{tabular}
  }
  }
  \end{center}
\caption{Exact meaning of the loss and error quantities reported in Section \ref{sec:numerics}.}\label{tab:reportedquantities}
\end{table}

\subsection{The parametric Klein-Gordon Equation}\label{SEC-KG}

We first test the Klein-Gordon equation parameterized by $(\alpha,\beta,\gamma)\in[-2,-1]\times[0,1]\times[0,1]$,
\begin{align}
\label{eq:kg}
\begin{split}
    u_{tt} + \alpha u_{xx} + \beta u + \gamma u^2 + x\cos{(t)} - x^2\cos^2{(t)} & = 0, \quad (x,t)\in[-1,1]\times[0,5],\\
    u(-1,t) = -\cos{(t)}, & \quad u(1,t)=\cos{(t)},\\
    u(x,0) & = x, \\
    u_t(x,0) & =0.
\end{split}
\end{align}
The full PINN is a $[2, 40, 40, 1]$-fully connected network with activation function $\cos{(z)}$ that is trained using uniformly distributed collocation points with $|\cC_o| = 10,000$, $|\cC_\partial| = 512$, $|\cC_i| = 512$. A learning rate of $0.0005$ is used with the ADAM optimizer and the maximum number of epochs being $75,000$. The parameter training set $\Xi_{\rm train}$ is a tensorial grid of size $10 \times 10 \times 10$ for a total of $1000$ parameter values. Up to $15$ neurons are generated by the greedy algorithm producing GPT-PINNs of sizes $[2, 1, 1]$ to $[2, 15, 1]$. The GPT-PINNs are trained at the same set of collocation points as the full PINN (i.e. $\cC_{pos}^r = \cC_{pos}$ for $pos \in \{o, \partial, i\}$) but with a learning rate of $0.025$ and (much smaller) $2000$ epochs. 
\begin{figure}[!htbp]
    \centering
    \includegraphics[width=0.32\textwidth]{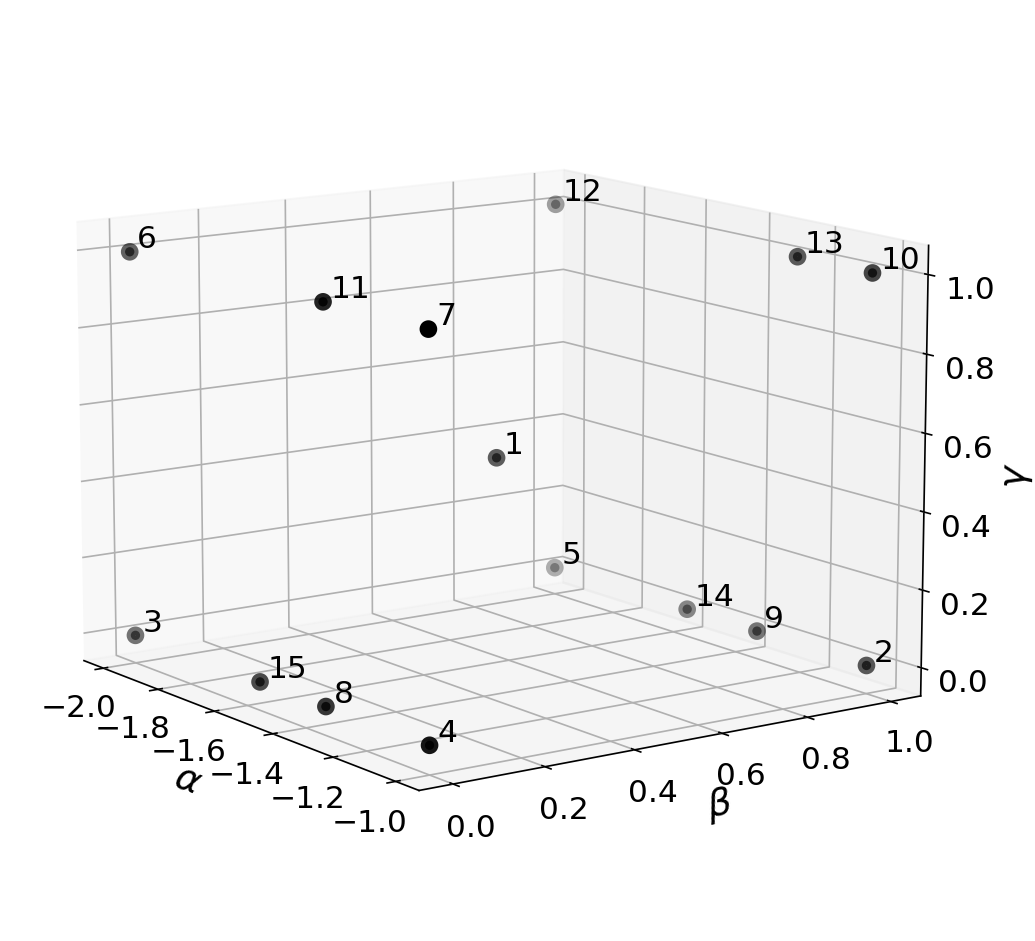}
    \includegraphics[width=0.32\textwidth]{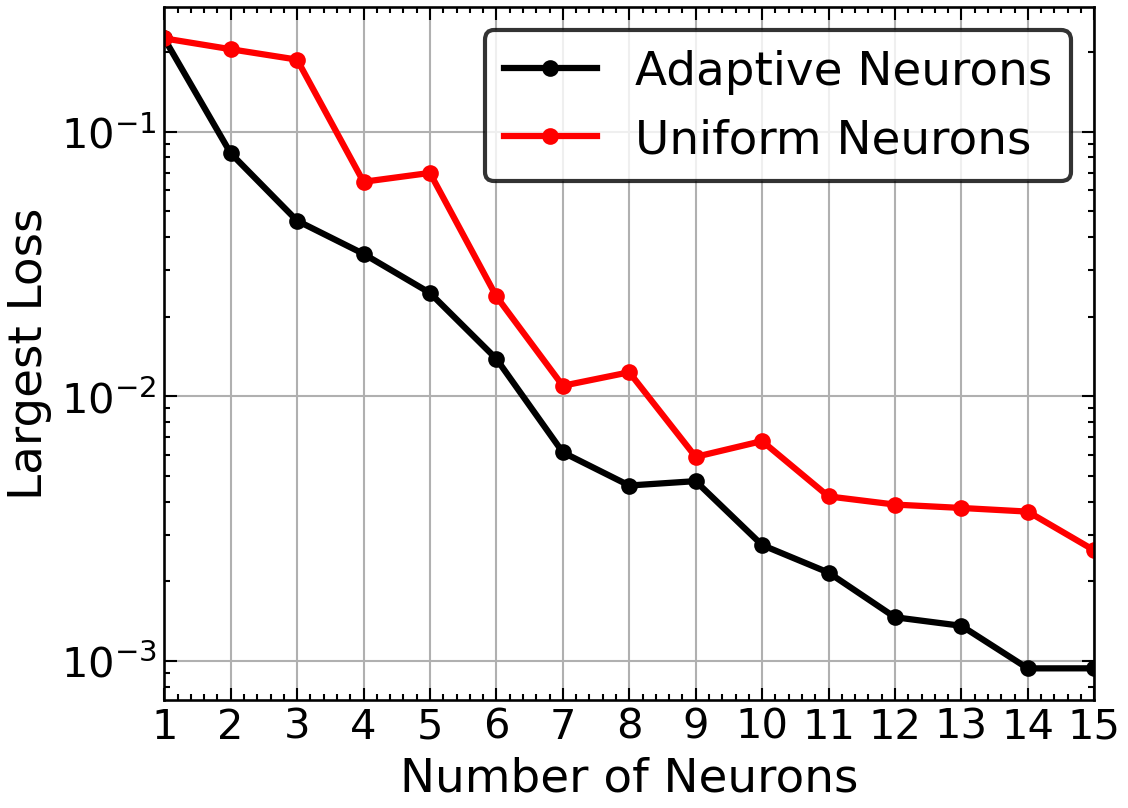}
    \includegraphics[width=0.32\textwidth]{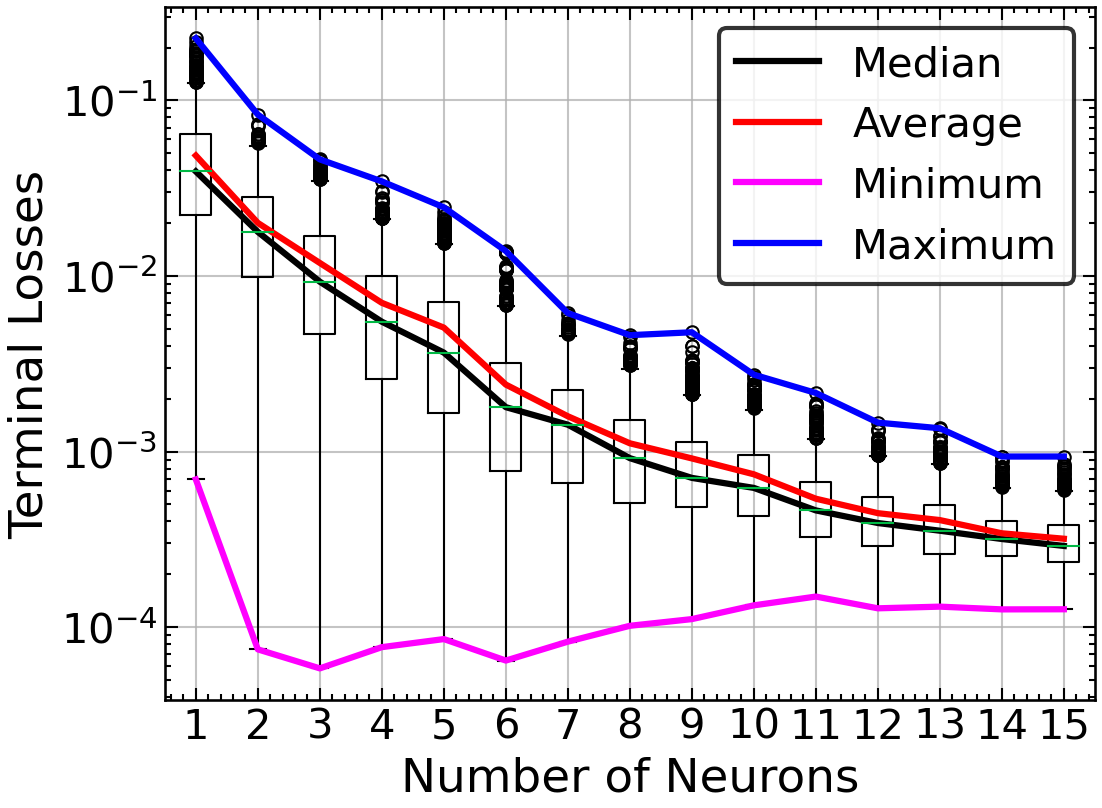}
    \caption{Klein-Gordon Equation training: The adaptively chosen parameter values (Left), worst-case GPT-PINN training losses (Middle), and the Box and Whisker plot of all adaptive GPT-PINN training losses (Right) during the outer-layer greedy training.}
    \label{fig:pinnloss_KG}
\end{figure}
\begin{figure}[!htbp]
    \centering
    \includegraphics[width=0.32\textwidth]{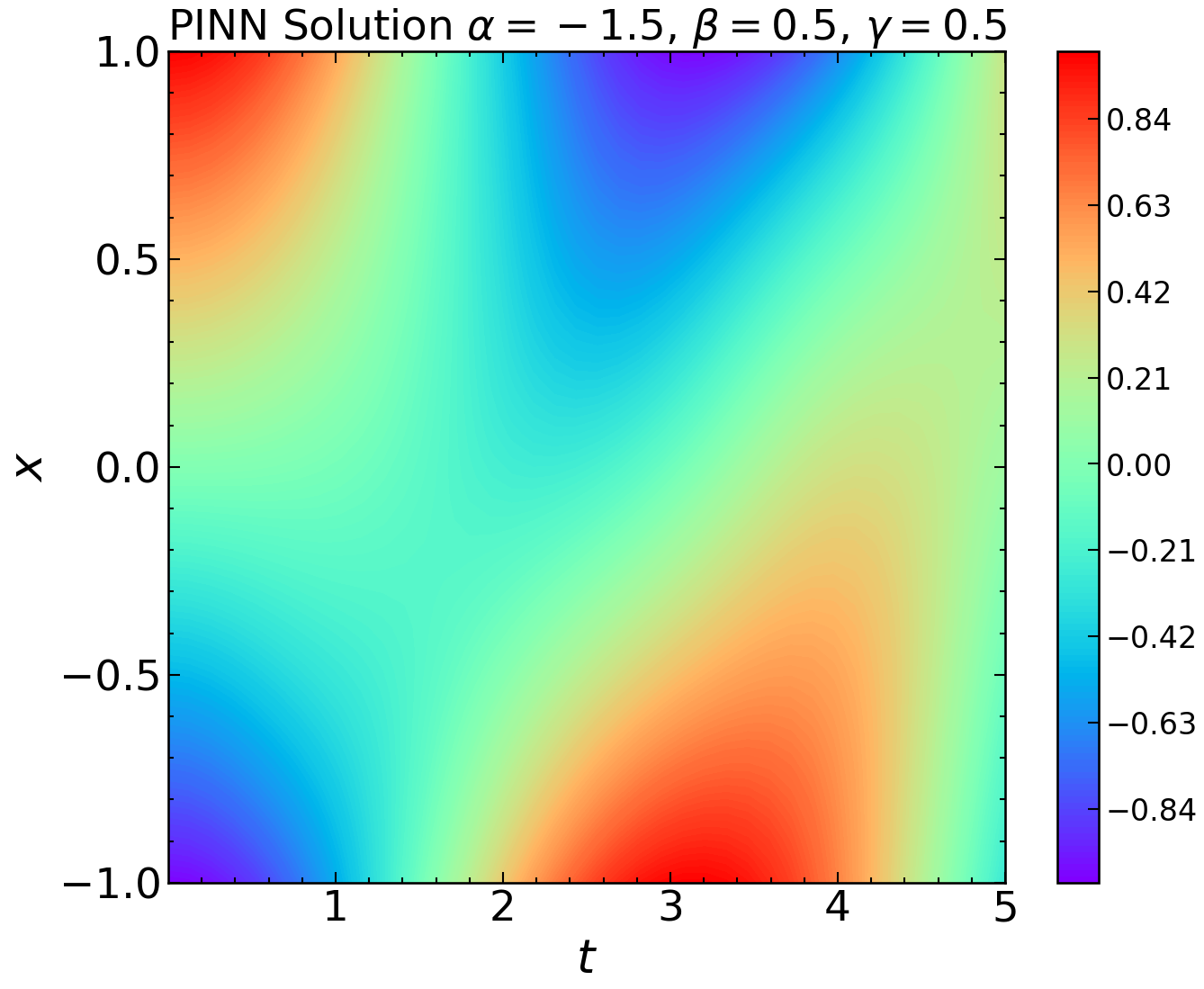}
    \includegraphics[width=0.32\textwidth]{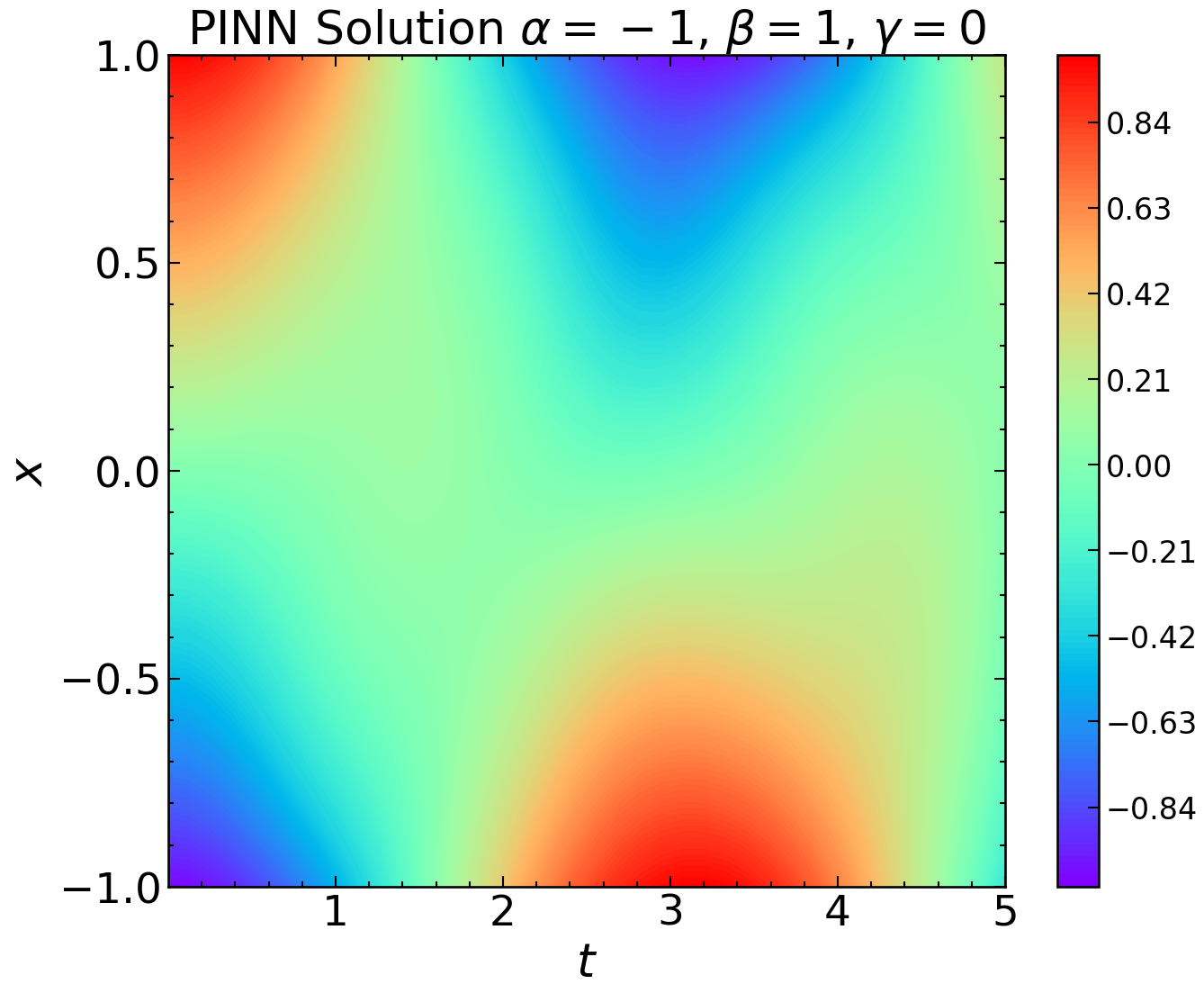}
    \includegraphics[width=0.32\textwidth]{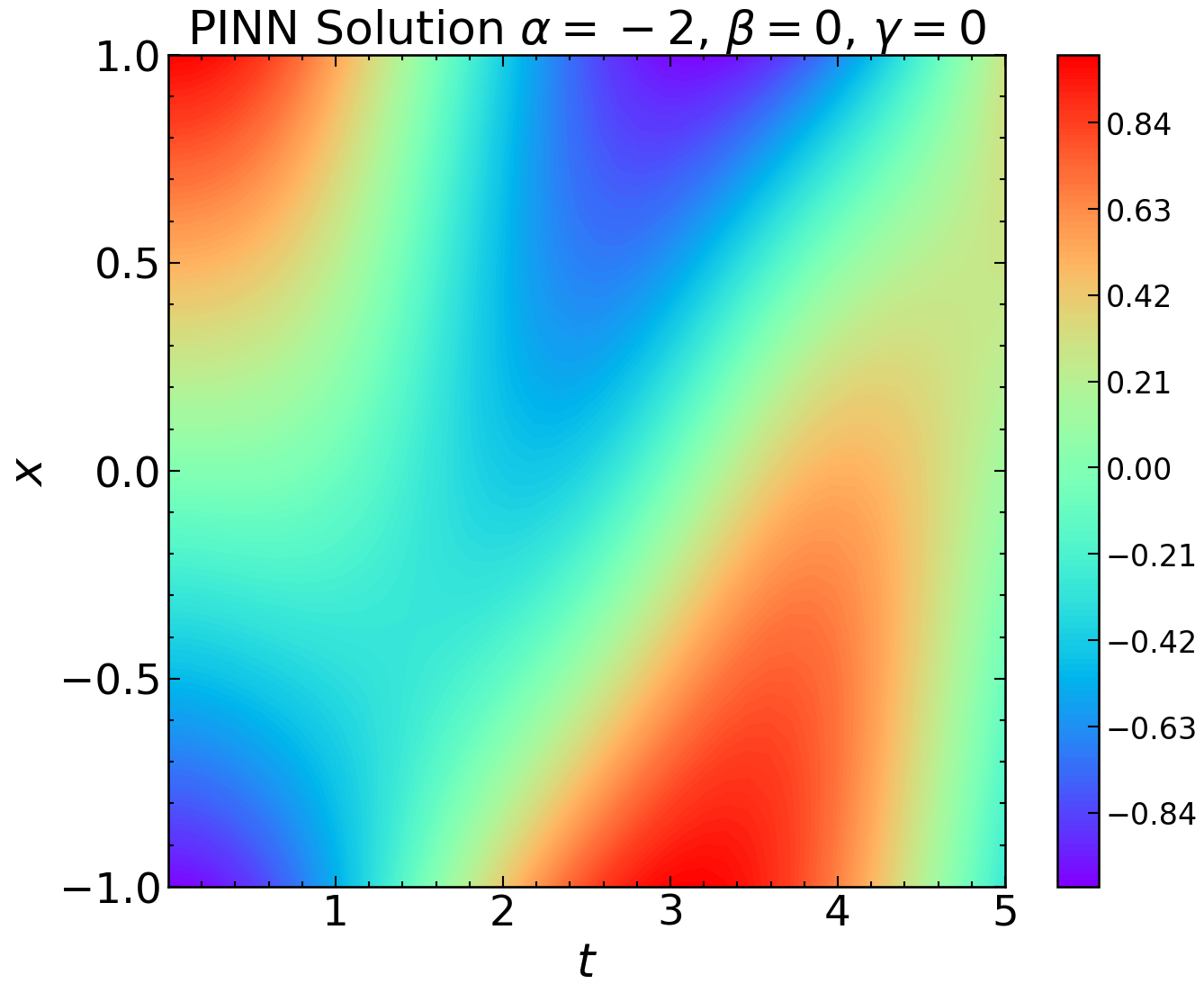}
    \caption{Klein-Gordon Equation: First three full PINN solutions found by the GPT-PINN that are used as the activation functions.}
    \label{fig:fullpinn_kg_sol}
\end{figure}

The GPT-PINN generates $15$ neurons, i.e. full PINNs at $\{(\alpha_i, \beta_i, \gamma_i)\}_{i=1}^{15}$. These parameter values and the worse-case offline training loss $\mathcal{L}_{\text{PINN}}^{\text{GPT}}(\bc(\bmu))$ after 2000 epochs as we increase the number of neurons (i.e. size of $\bc(\bmu)$) in the hidden layer of GPT-PINN are shown in Figure \ref{fig:pinnloss_KG}. Figure \ref{fig:fullpinn_kg_sol} shows the first three PINN solutions adaptively selected by GPT-PINN. 
 It is clear that the sampled parameter values are toward the boundaries of the domain and that the decrease in training loss is exponential. Both features are consistent with typical RBM results. Moreover, we emphasize that to achieve 3 digits of accuracy across the parameter domain, we only need to train the full PINN 15 times. This contrast with pure data-driven approaches inherits that of RBM with POD approaches in that RBM requires much less full-order solves. In comparison, we sample the parameter domain uniformly (i.e. without the greedy approach of GPT-PINN), it is clear from Figure \ref{fig:pinnloss_KG} Middle that the adaptive ``learned neurons'' performs 2 to 3 times better than the non-adaptive ``uniform neurons''. The fact that the latter performs reasonably well underscores the power of our novel idea of using pre-trained PINNs as activation functions.

\begin{figure}[!htbp]
    \centering
    \includegraphics[width=0.32\textwidth]{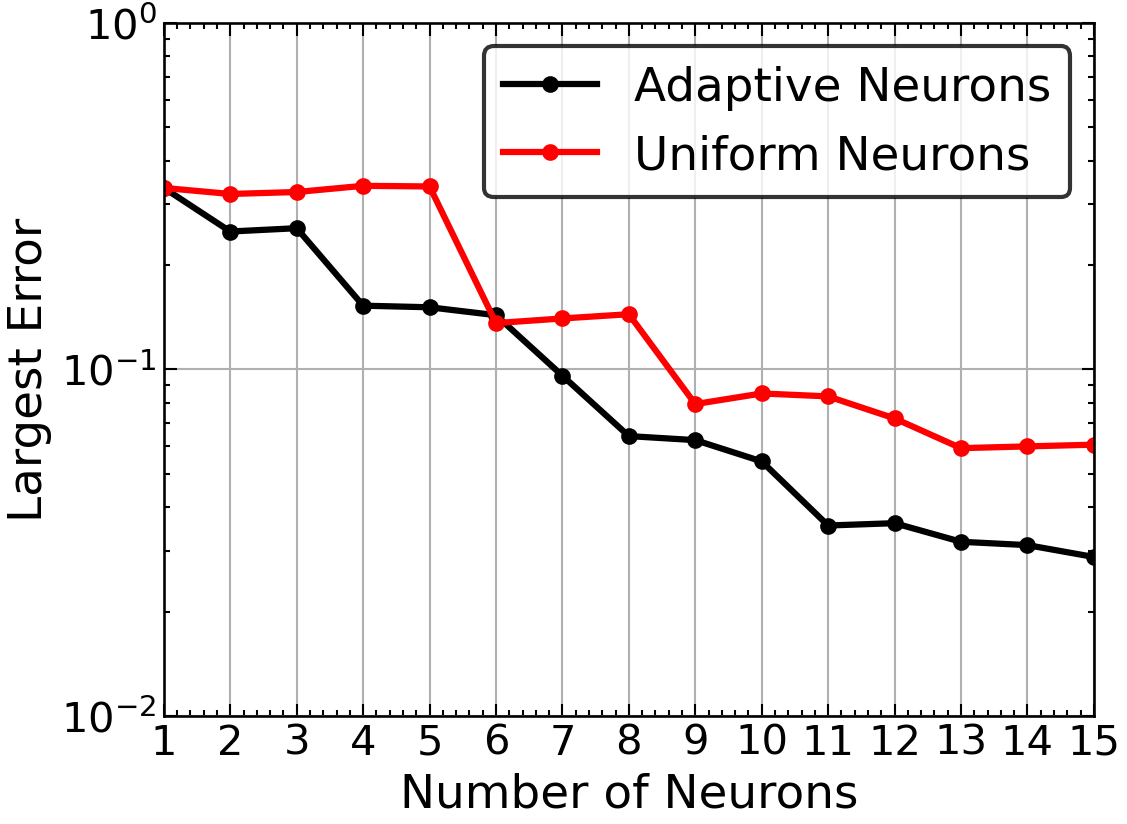}
    \includegraphics[width=0.32\textwidth]{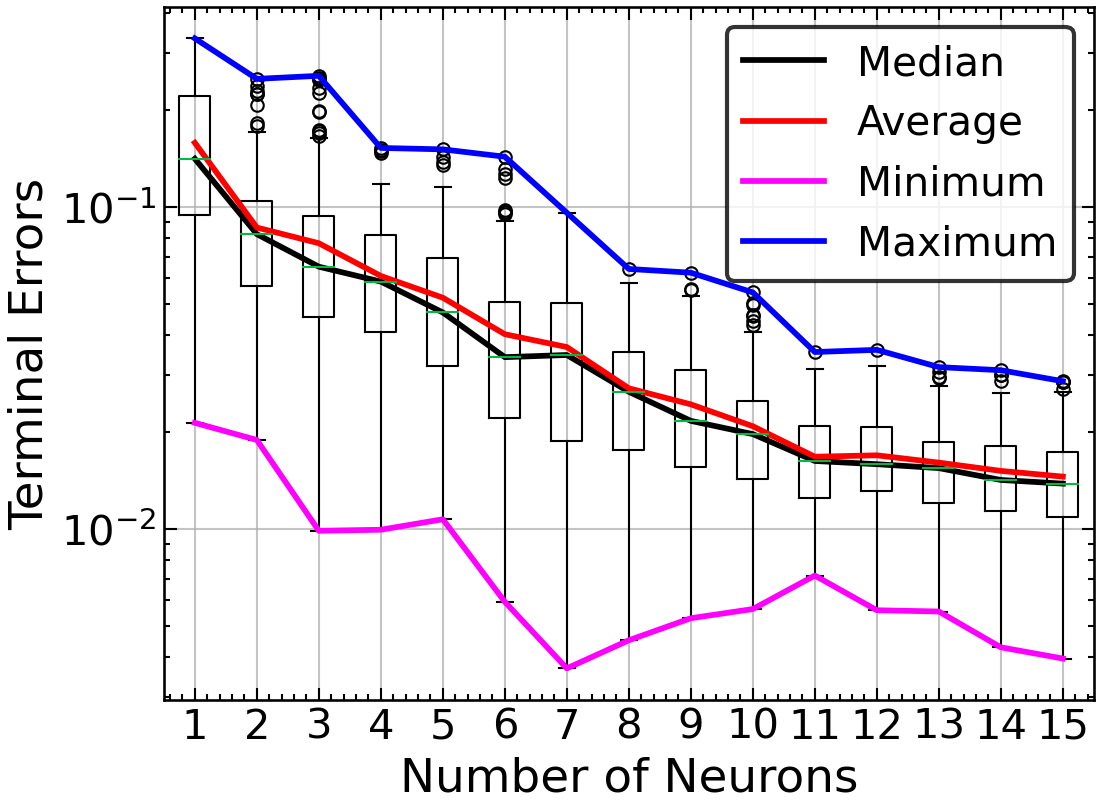}
    \includegraphics[width=0.32\textwidth]{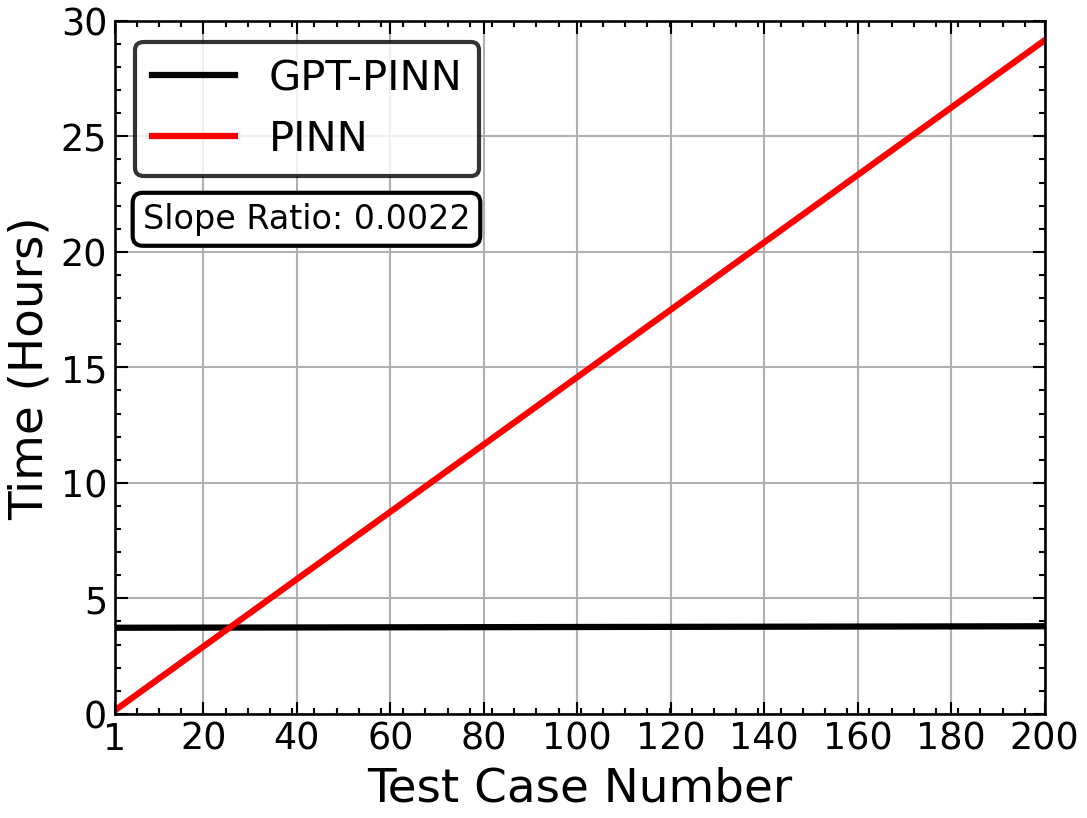}
    \caption{Klein-Gordon Equation testing: Worst-case test error of the GPT-PINN of various sizes (Left), Box and Whisker plot of all adaptive GPT-PINN testing errors (Middle), and cumulative run time of the full PINN versus the GPT-PINN (Right).}
    \label{fig:pinnerror_KG}
\end{figure}

Next, we test the GPT-PINN on $\Xi_{\rm test}$ consisting of $200$ randomly selected parameter values distinct from the adaptively chosen ``learned neurons'' and ``uniform neurons''. Figure \ref{fig:pinnerror_KG} displays the largest error for each size of the GPT-PINN. The trend is again exponential. Finally, to show the efficiency of the method, we plot in Figure \ref{fig:pinnerror_KG} Right the cumulative run-time when both the full PINN and the (reduced) GPT-PINN are repeatedly called. The starting point of the GPT-PINN line reflects all offline preparation time. It is clear that the GPT-PINN line increases very slowly reflecting the fact that its marginal cost is tiny. In fact, it is about 0.0022 of that of the full PINN. The intersection points reflect how many simulations would it be worthwhile to invest in GPT-PINN.  We remark that future work includes driving the intersection point down to essentially comparable to the number of neurons in GPT-PINN, which is the absolute minimum it could be.

\begin{figure}[!htbp]
    \centering
    \includegraphics[width=0.32\textwidth]{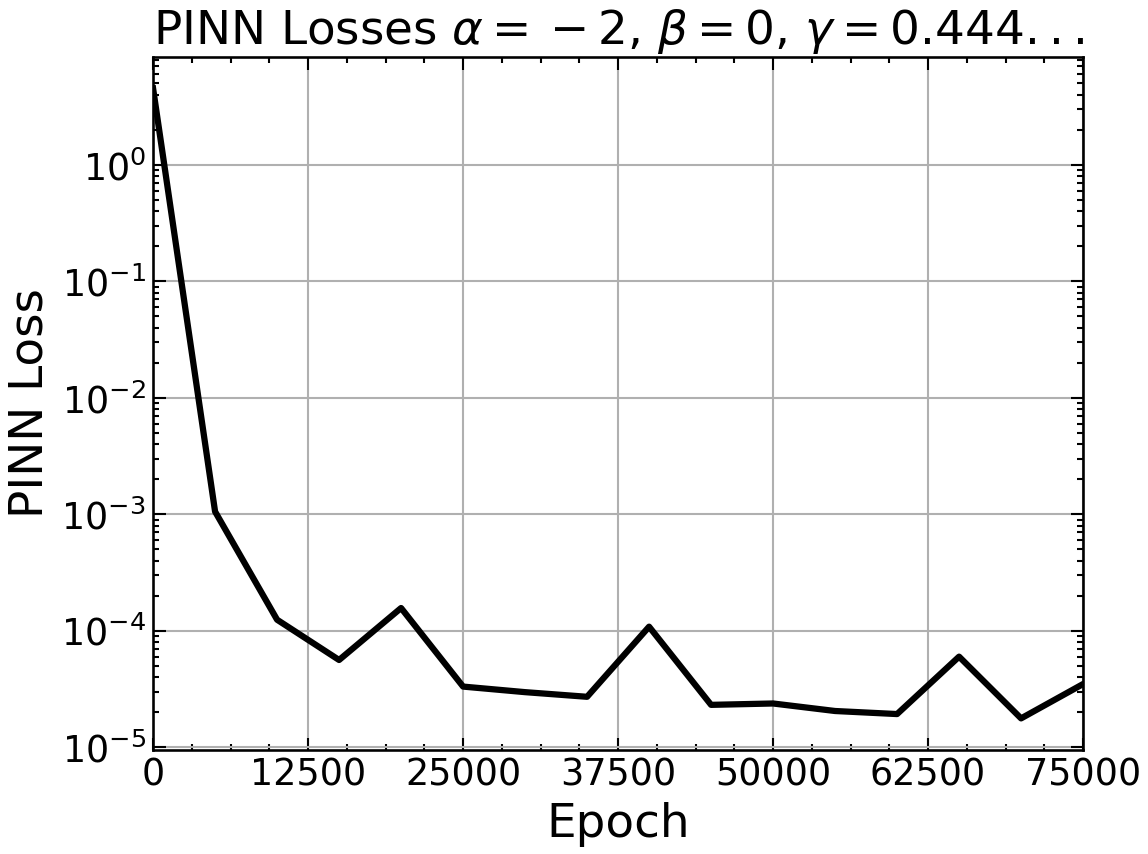}
    \includegraphics[width=0.32\textwidth]{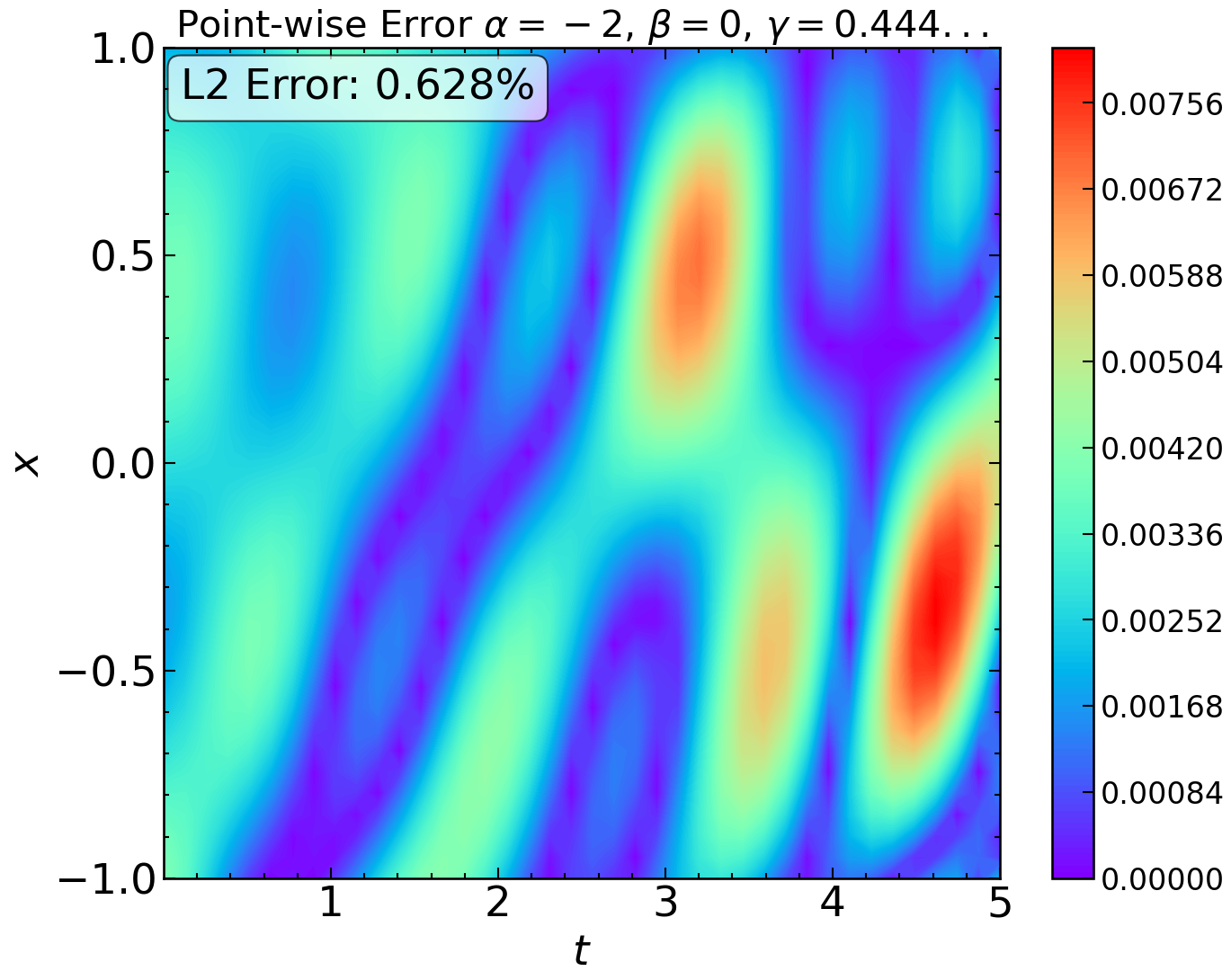}
    \includegraphics[width=0.32\textwidth]{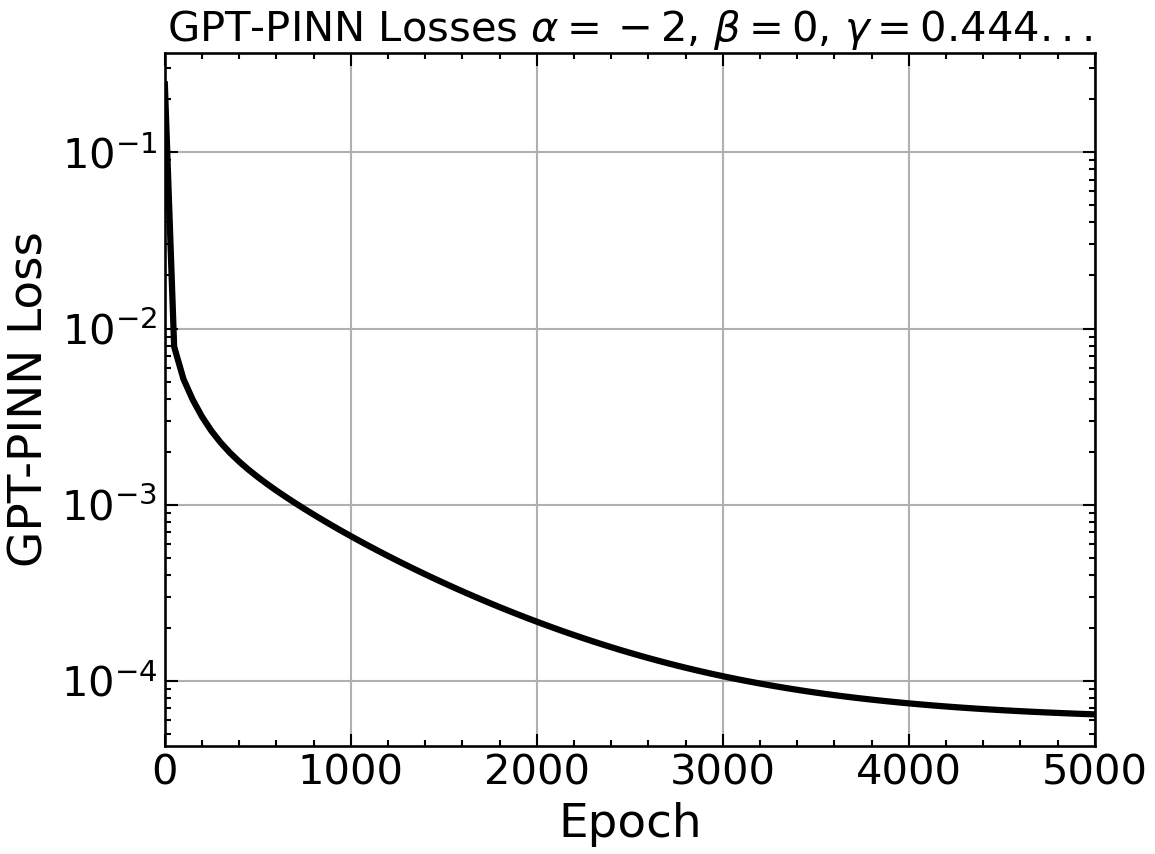}\\
    \includegraphics[width=0.32\textwidth]{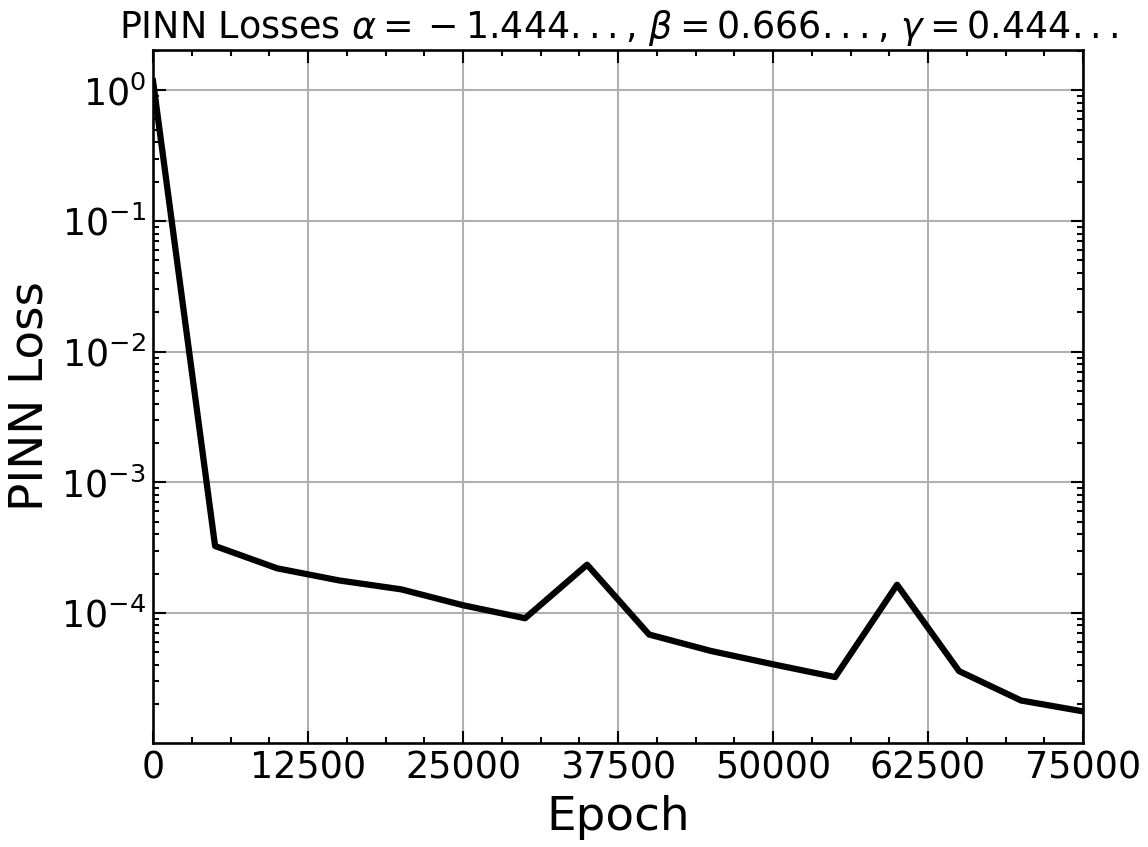}
    \includegraphics[width=0.32\textwidth]{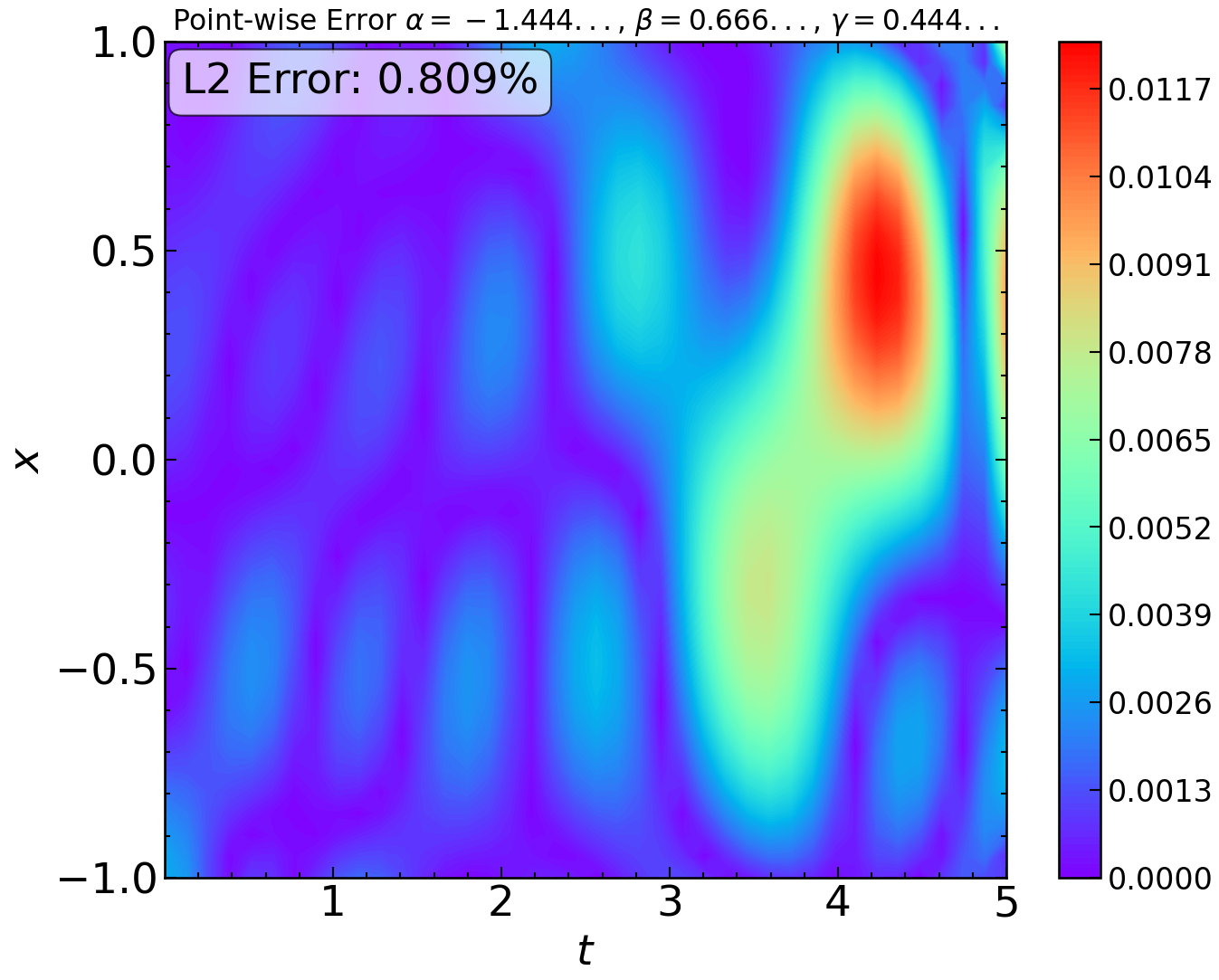}
    \includegraphics[width=0.32\textwidth]{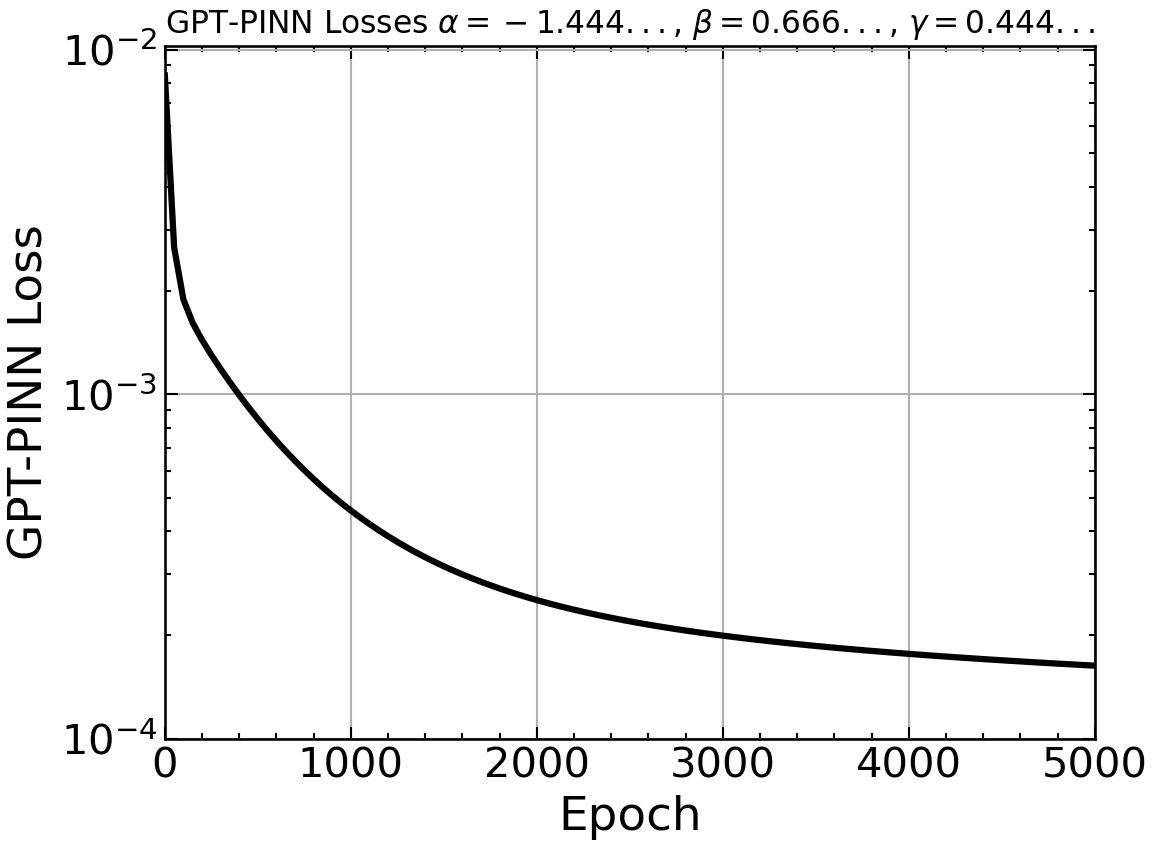}
    \caption{Klein-Gordon Equation: Full PINN training loss (Left) and GPT-PINN training loss (Right) as functions of epochs for various parameters. Plotted in the middle are the point-wise errors of the corresponding GPT-PINN solution.}
    \label{fig:loss_v_epoch_kg}
\end{figure}

Last but not least, we show the training losses as functions of epochs in Figure \ref{fig:loss_v_epoch_kg} for both the full PINN and GPT-PINN. We note the interesting phenomenon that the GPT-PINN loss decreases more smoothly than the full PINN. To give a sense of the error distribution, we also plot the point-wise error of the GPT-PINN solution.

\subsection{The parametric viscous Burgers' Equation}\label{SEC-B}

Next, we test GPT-PINN on the Burgers' equation with one parameter, the viscosity $\nu\in[0.005, 1]$.
\begin{align}
\label{eq:burgers}
\begin{split}
    u_t + uu_x - \nu u_{xx} & = 0, \quad (x,t)\in[-1,1]\times[0,1],\\
    u(-1,t)=u(1,t) & =0,\\
    u(x,0) & = -\sin{(\pi x)}.
\end{split}
\end{align}

The full PINN is a $[2, 20, 20, 20, 20, 1]$-fully connected network with activation function $\tanh(z)$ that is trained  using uniformly distributed collocation points with $|\cC_o| = 10,000$, $|\cC_\partial| = 100$, $|\cC_i| = 100$. A learning rate of $0.005$ is used with the ADAM optimizer. A maximum number of $60,000$ epochs is run with a stopping criteria of $2\times10^{-5}$ implemented on the loss values. The parameter training set is a uniform grid of size $129$ in the $\nu$-domain. Up to $9$ neurons are generated by the greedy algorithm producing the reduced GPT-PINNs of sizes $[2, 1, 1]$ to $[2, 9, 1]$. The GPT-PINNs are trained at the same set of collocation points as the full PINN but with a learning rate of $0.02$ and $2000$ epochs. The solutions of \cref{eq:burgers} develop near-discontinuities as time evolves when $\nu$ is small. In this scenario,  $\left(\Psi^{\theta^i}_{\mathsf{NN}}\right)_{xx}$ is of little value in the training of GPT-PINN when $x$ is close to these large gradients. We therefore exclude the collocation points where  $\left|\left(\Psi^{\theta^i}_{\mathsf{NN}}\right)_{xx}\right|$ is within the top $20\%$ of all such values. That is
\[
\cC_{pos}^r = \cC_{pos} \backslash \left\{x: \left|\left(\Psi^{\theta^i}_{\mathsf{NN}}\right)_{xx}(x)\right| > 0.8 \max_x \left|\left(\Psi^{\theta^i}_{\mathsf{NN}}\right)_{xx}(x)\right|\right\}, \quad pos \in \{o, \partial, i\}.
\]

\begin{figure}[!htbp]
    \centering
    \includegraphics[width=1\textwidth]{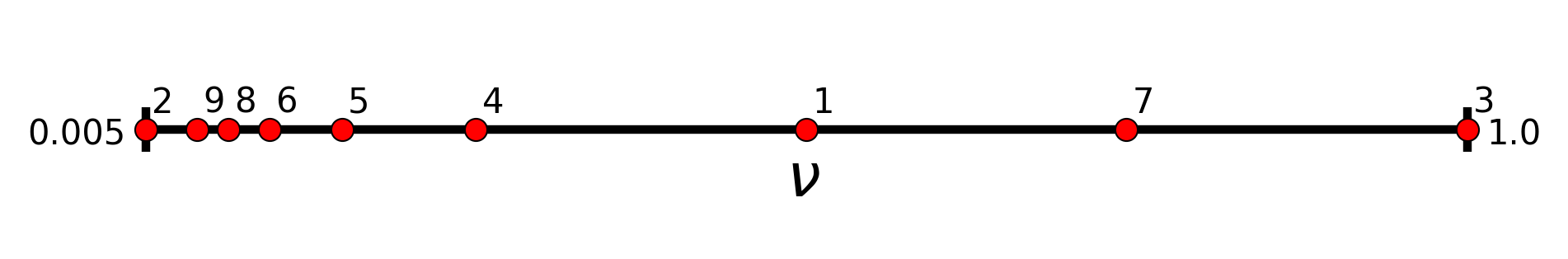}
    \includegraphics[width=0.49\textwidth]{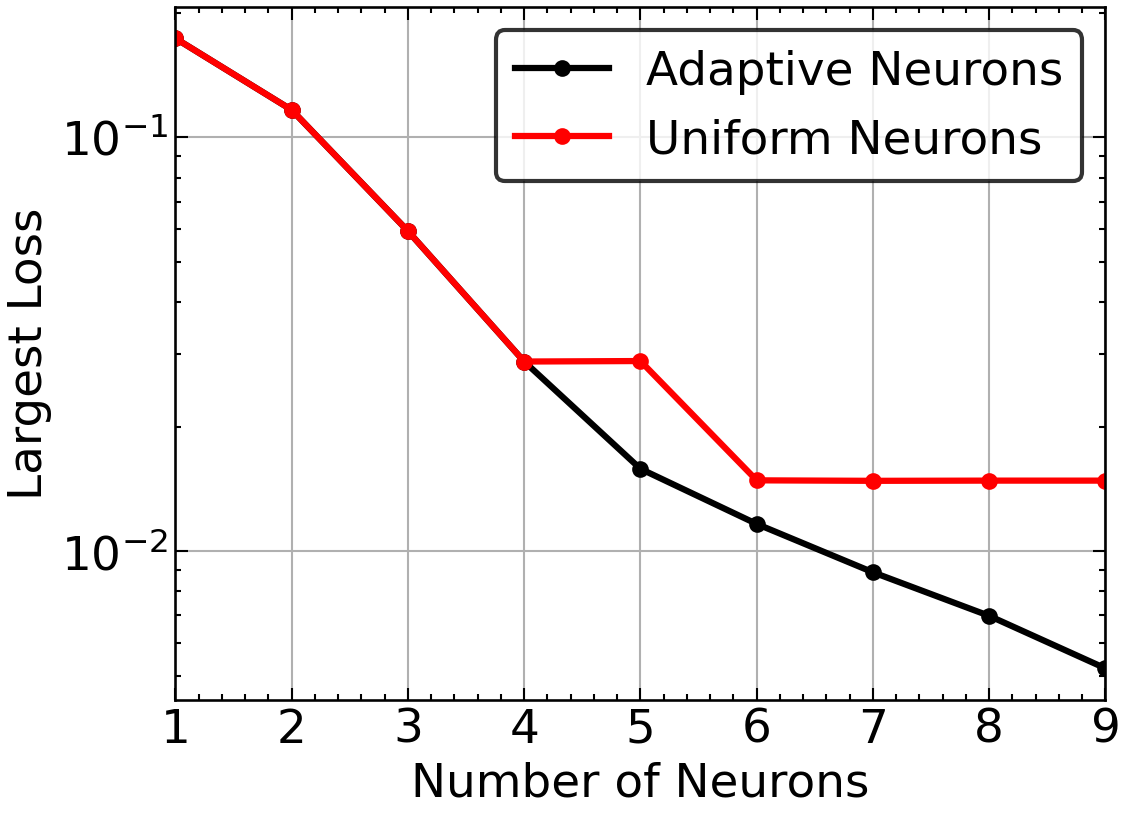}
    \includegraphics[width=0.49\textwidth]{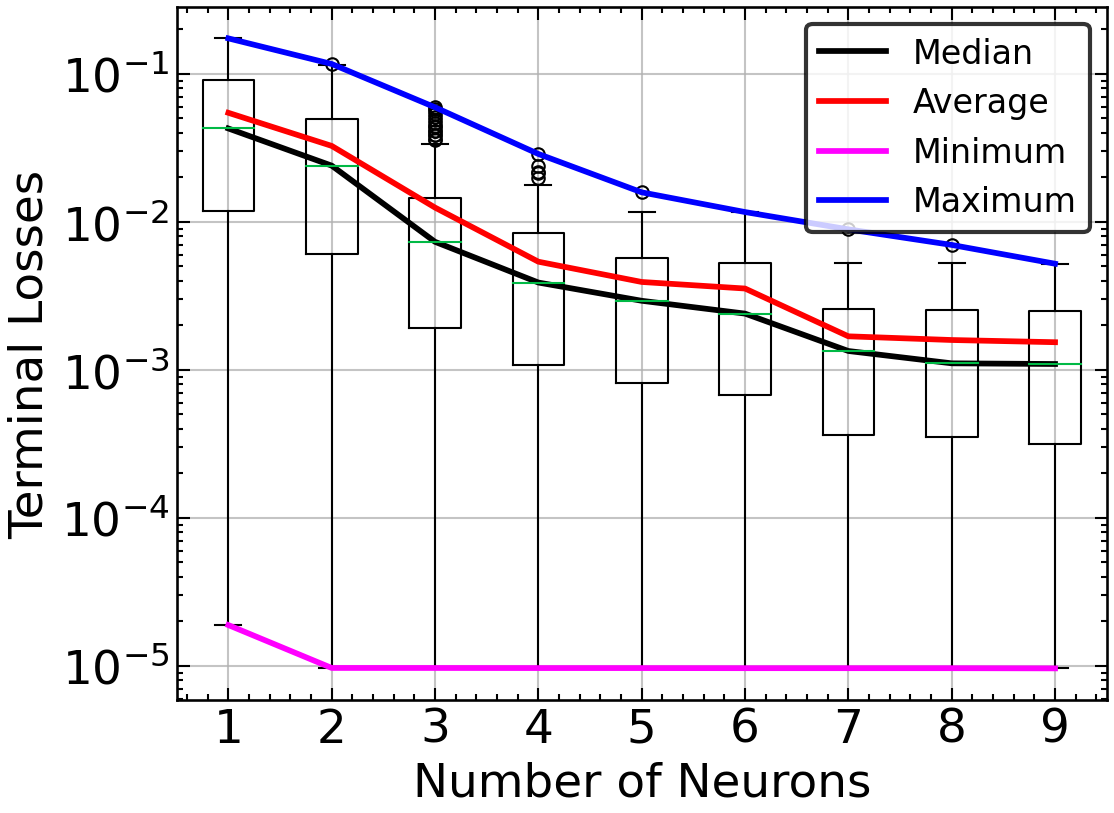}
    \caption{Burgers' Equation training: The adaptively chosen parameter values (Top), worst-case GPT-PINN training losses (Bottom Left), and the Box and Whisker plot of all GPT-PINN training losses (Bottom Right) during the outer-layer greedy training.}
    \label{fig:pinnloss_B}
\end{figure}
\begin{figure}[!htbp]
    \centering
    \includegraphics[width=0.32\textwidth]{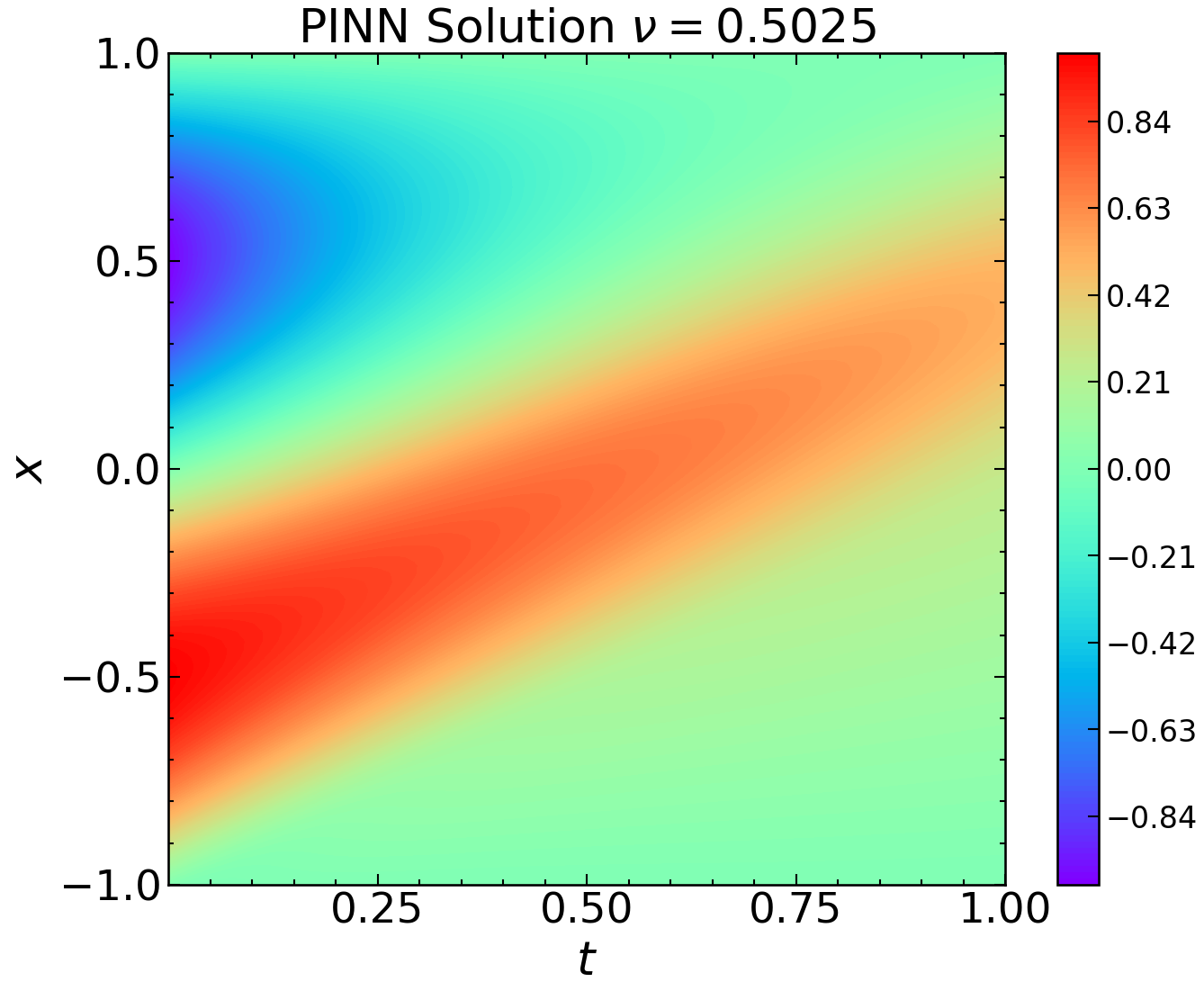}
    \includegraphics[width=0.32\textwidth]{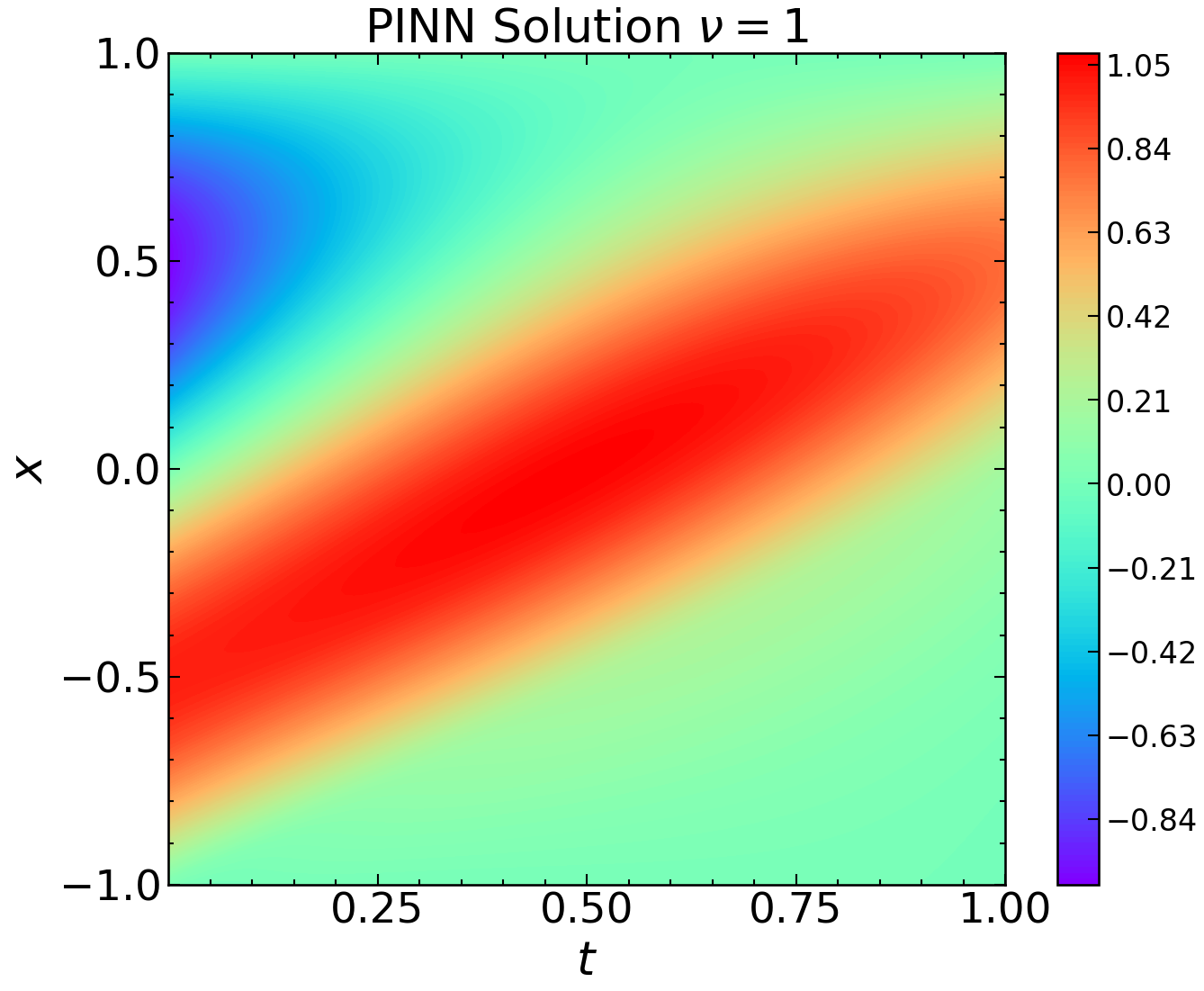}
    \includegraphics[width=0.32\textwidth]{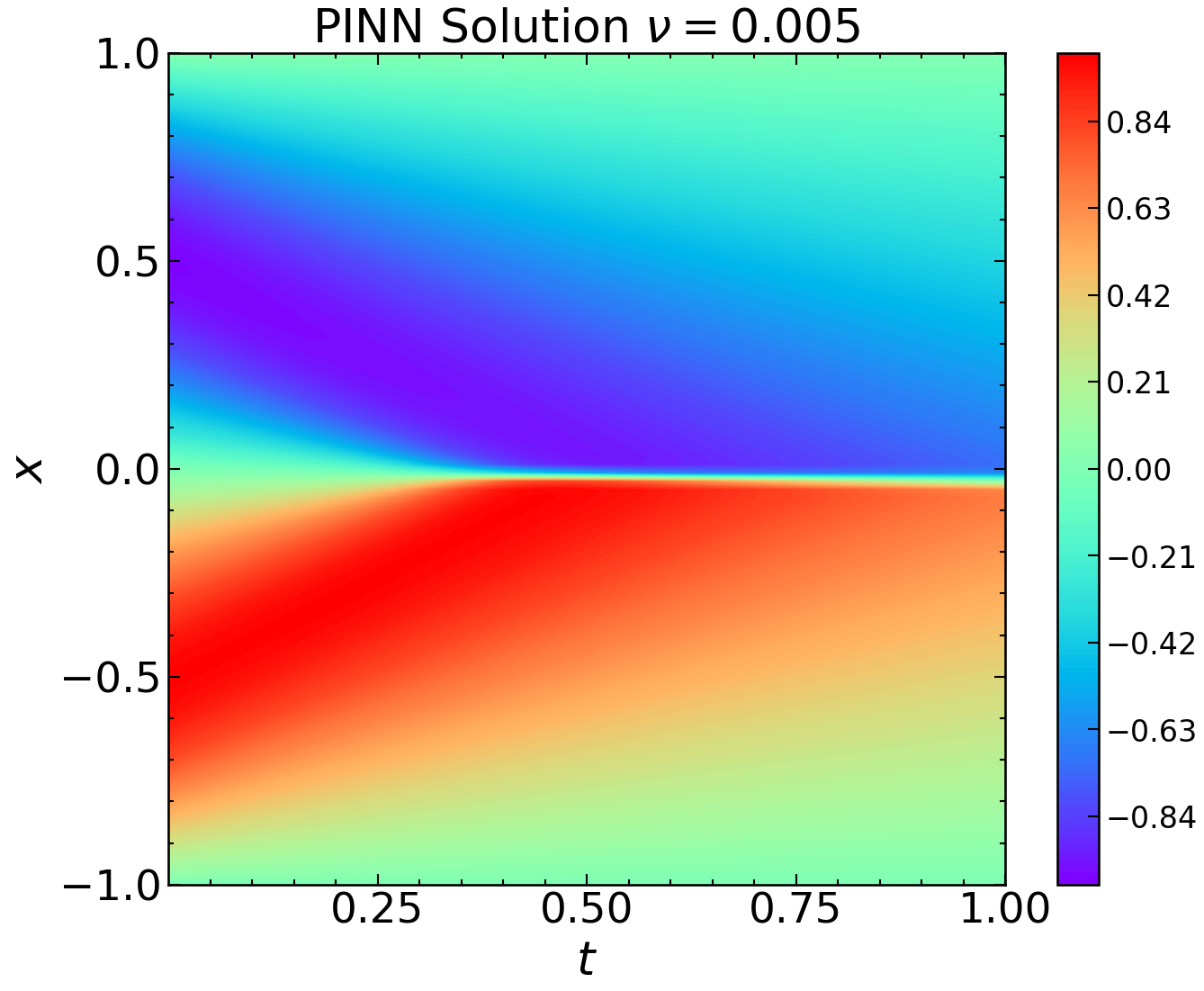}
    \caption{Burgers' Equation: First three full PINN solutions found by the GPT-PINN that are used as the activation functions.}
    \label{fig:fullpinn_b_sol}
\end{figure}

The GPT-PINN generates $9$ neurons, i.e. full PINNs at $\{(\nu_i, \}_{i=1}^{9}$. These parameter values and the worse-case offline training loss $\mathcal{L}_{\text{PINN}}^{\text{GPT}}(\bc(\bmu))$ after 2000 epochs as we increase the number of neurons (i.e. size of $\bc(\bmu)$) in the hidden layer of GPT-PINN are shown in Figure \ref{fig:pinnloss_B}. Figure \ref{fig:fullpinn_b_sol} shows the first three PINN solutions adaptively selected by GPT-PINN. We observe behavior that is similar to the Klein-Gordon case and consistent with typical RBM results. The adaptive ``learned neurons'' again perform 3 to 4 times better than the non-adaptive ``uniform neurons'' which already perform reasonably well, underscoring the power of our novel idea of using pre-trained PINNs as activation functions.

\begin{figure}[!htbp]
    \centering
    \includegraphics[width=0.32\textwidth]{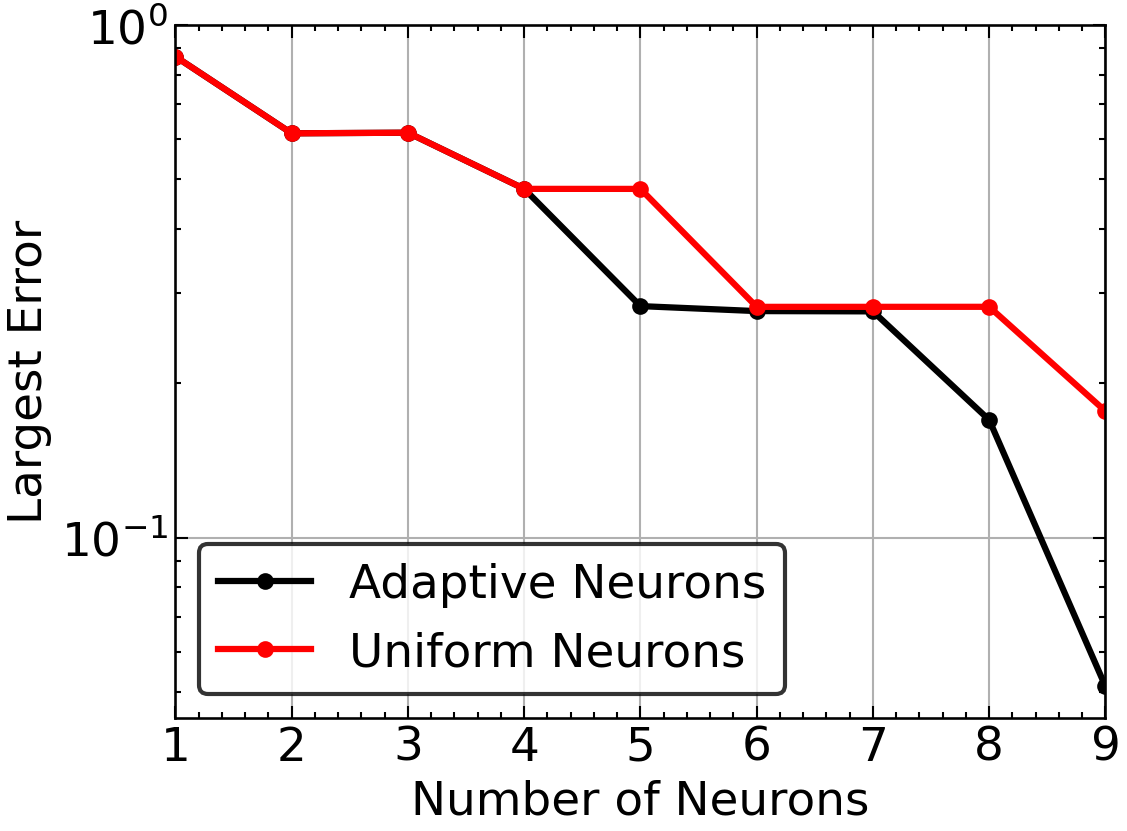}
    \includegraphics[width=0.32\textwidth]{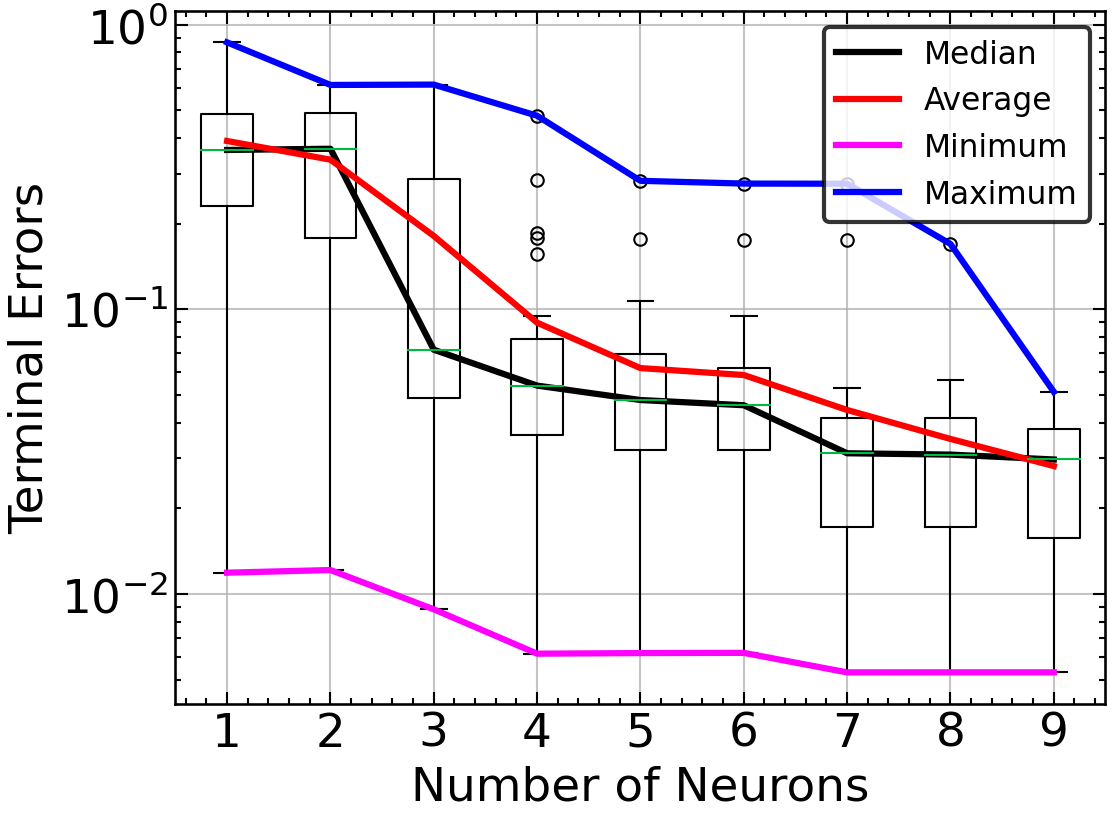}
    \includegraphics[width=0.32\textwidth]{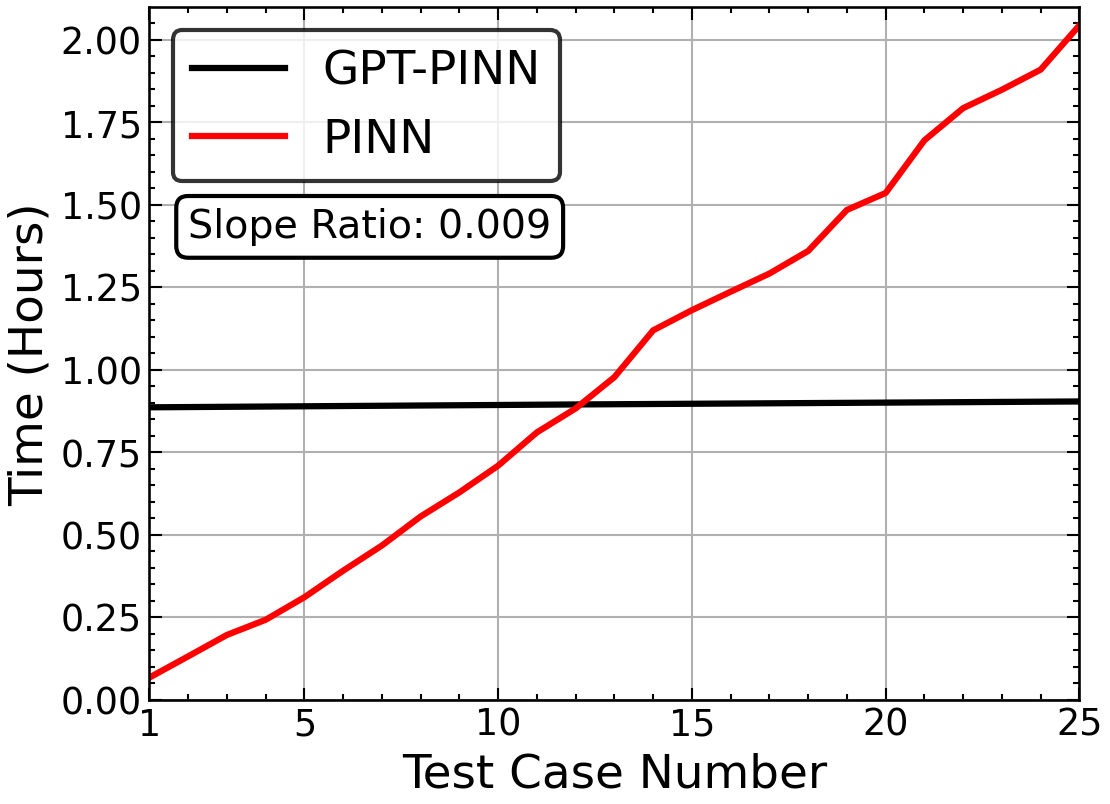}
    \caption{Burgers' Equation testing: Worst-case test error of the GPT-PINN of various sizes (Left), Box and Whisker plot of all adaptive GPT-PINN testing errors (Middle), and cumulative run time of the full PINN versus the GPT-PINN (Right).}
    \label{fig:pinnerror_B}
\end{figure}
\begin{figure}[!htbp]
    \centering
    \includegraphics[width=0.32\textwidth]{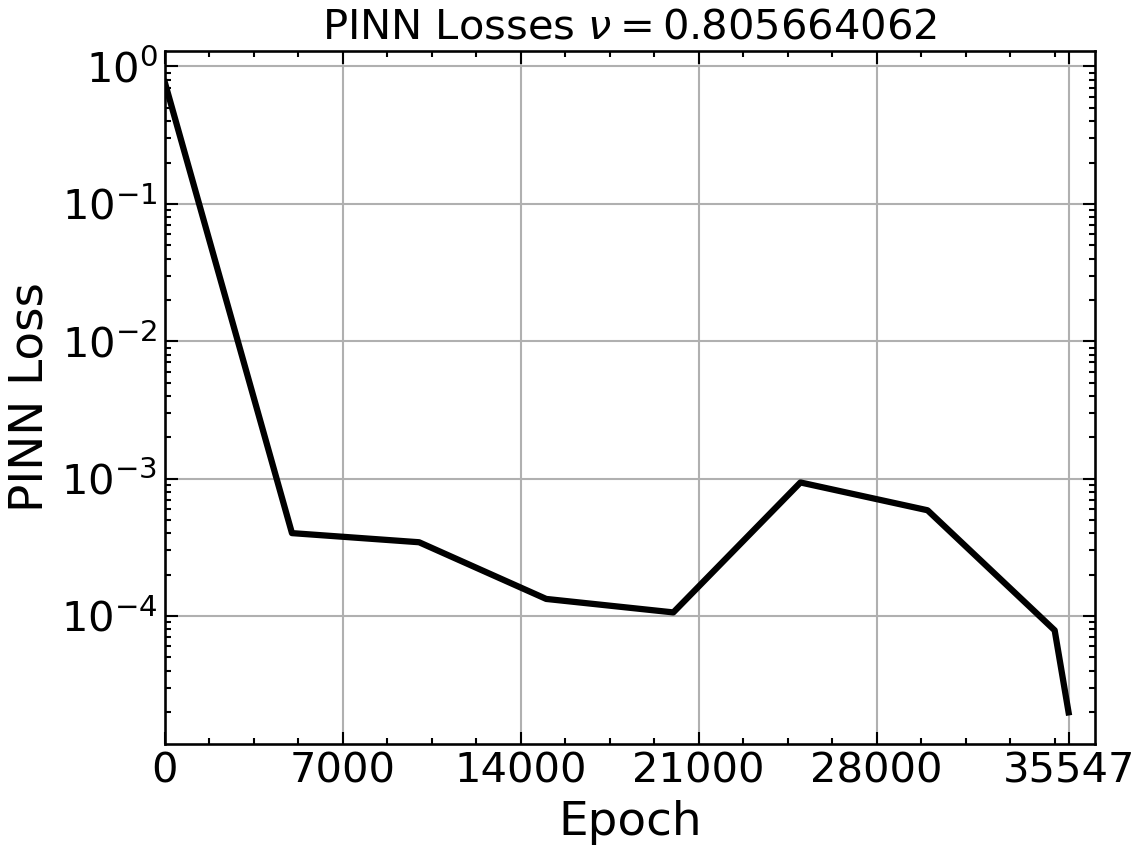}
    \includegraphics[width=0.32\textwidth]{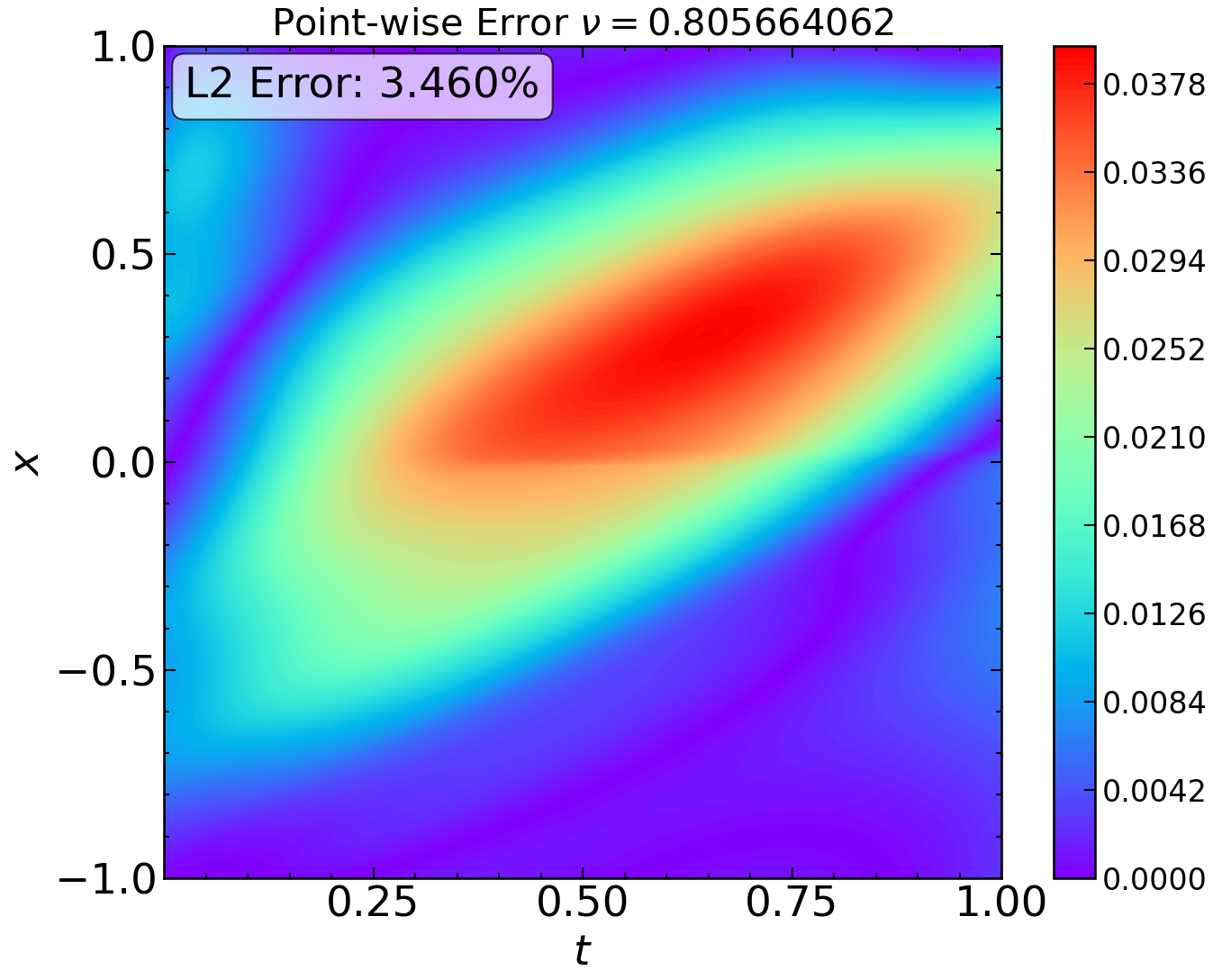}
    \includegraphics[width=0.32\textwidth]{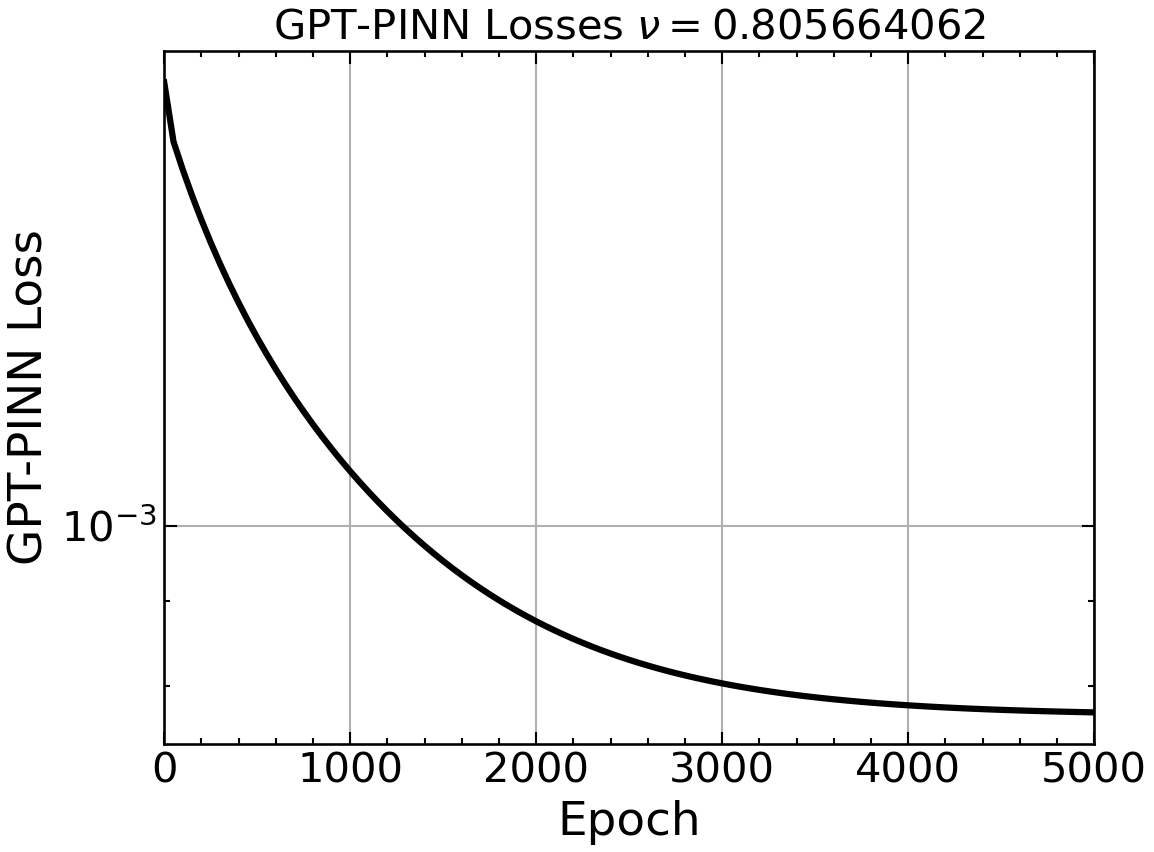}
    \includegraphics[width=0.32\textwidth]{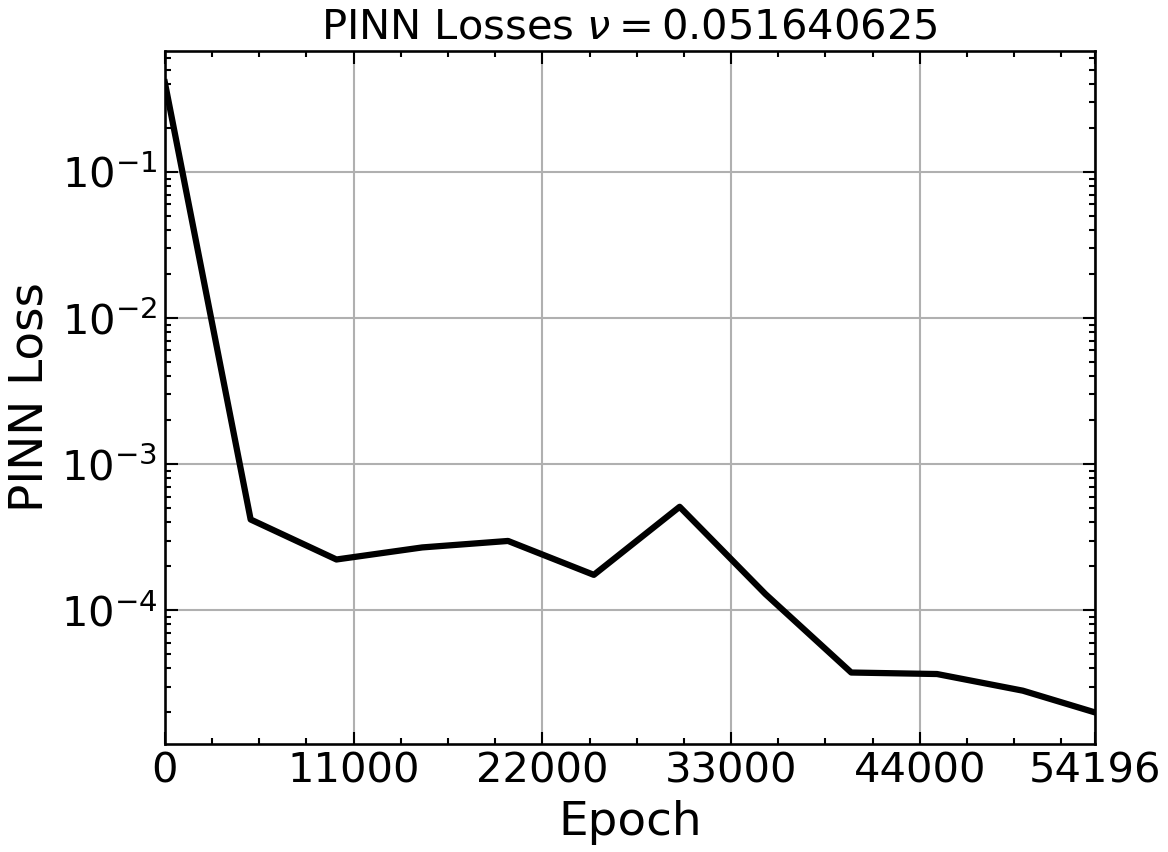}
    \includegraphics[width=0.32\textwidth]{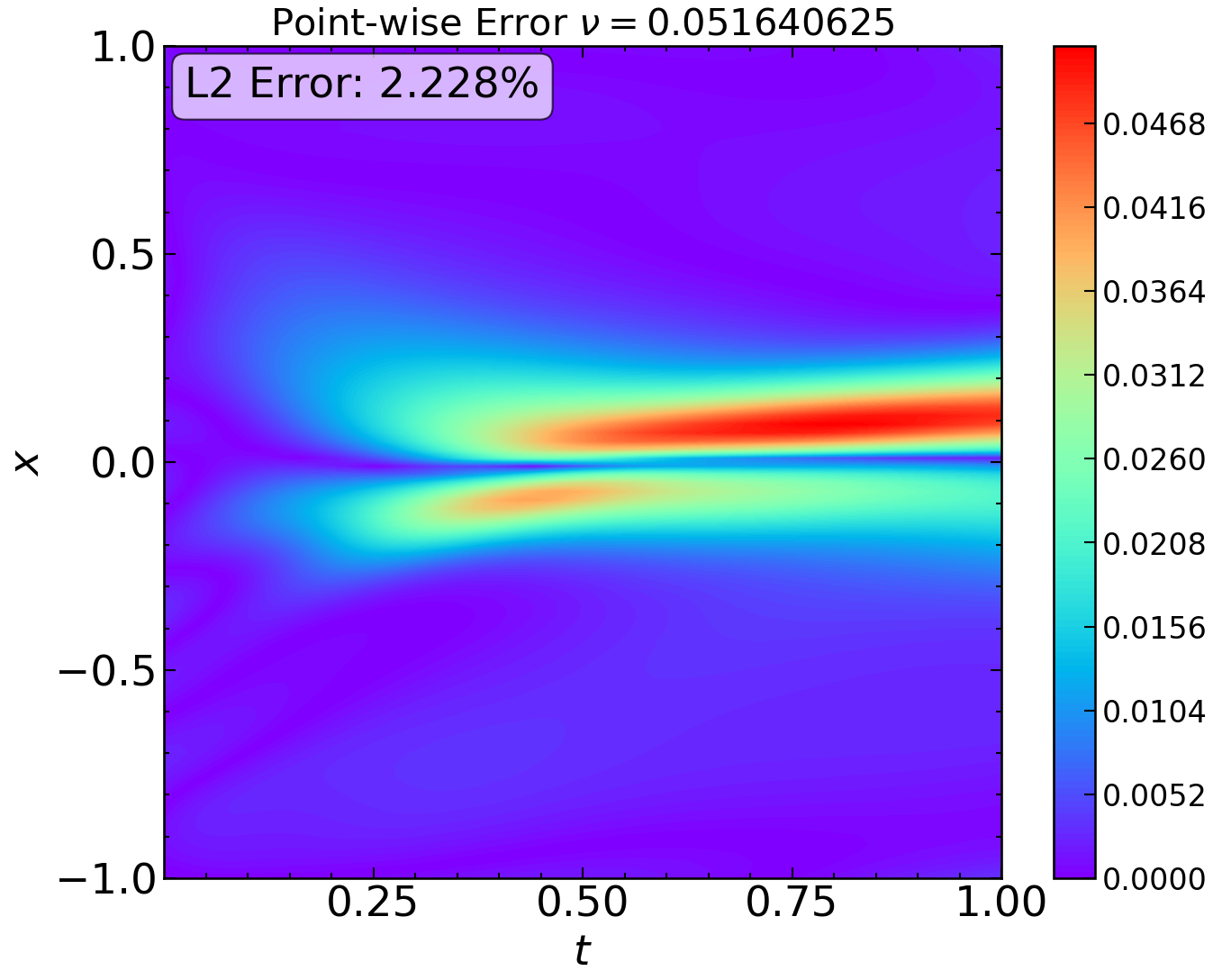}
    \includegraphics[width=0.32\textwidth]{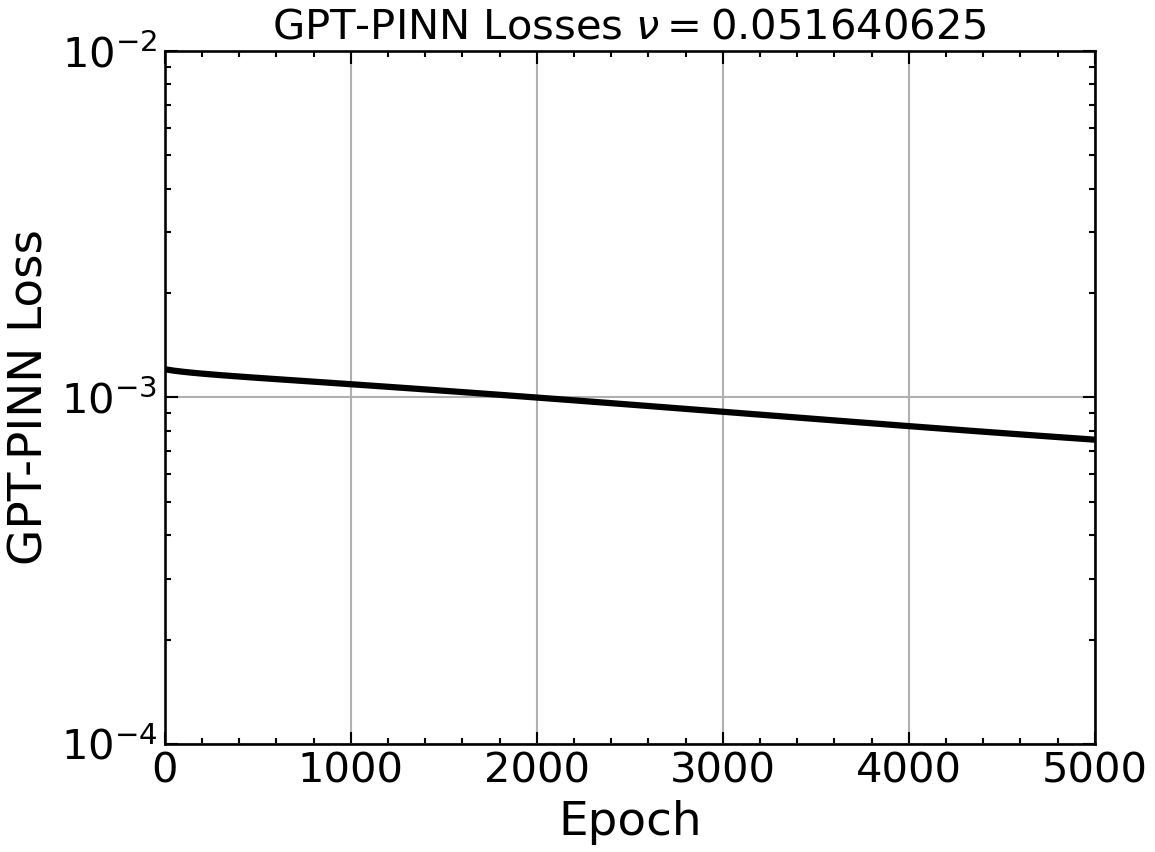}
    \caption{Burgers' Equation: Full PINN training loss (Left) and GPT-PINN training loss (Right) as functions of epochs for various parameters. Plotted in the middle are the point-wise errors of the corresponding GPT-PINN solution.}
    \label{fig:loss_v_epoch_B}
\end{figure}
Next, we test the GPT-PINN on $25$ parameter values. Figure \ref{fig:pinnerror_B} displays the largest error for each size of the GPT-PINN. The trend is again exponential. Finally, to show the efficiency of the method, we plot in Figure \ref{fig:pinnerror_B} Right the cumulative run-time when both the full PINN and the (reduced) GPT-PINN are repeatedly called. It is clear that the GPT-PINN line increases very slowly (a relative speed of $0.009$ in comparison to the full PINN) and that it is worthwhile to invest in GPT-PINN for a very modest number ($12$) of queries. 
We again show the training losses as functions of epochs in Figure \ref{fig:loss_v_epoch_B} for both the full PINN and GPT-PINN. We note again that the GPT-PINN loss decreases more smoothly than the full PINN. This result also verifies the efficacy of our initialization strategy since the starting loss of the GPT-PINN is already very low.

\subsection{The parametric Allen-Cahn Equation}\label{SEC-AC}

Finally, we test the Allen-Cahn equation parameterized by $(\lambda,\epsilon)\in[0.0001,0.001]\times[1,5]$
\begin{align}
\label{eq:AC}
\begin{split}
    u_t - \lambda u_{xx} + \epsilon (u^3-u) & =0, \quad (x,t)\in[-1,1]\times[0,1]\\
    u(-1,t)=u(1,t) & =-1\\
    u(x,0)&=x^2\cos{(\pi x)}.
\end{split}
\end{align}

The SA-PINN \cite{mcclenny2020self} is a $[2, 128, 128, 128, 128, 1]$-fully connected network with activation function $\tanh(z)$ that is trained on collocation points distributed by a Latin hypercube sampling with $|\cC_o| = 20,000$, $|\cC_\partial| = 100$, $|\cC_i| = 512$. A learning rate of $0.005$ with $10,000$ epochs of ADAM optimization followed by $10,000$ epochs of L-BFGS optimization with a learning rate of $0.8$ is used. The parameter training set is a grid of size $121$ uniform parameter values. Up to $9$ neurons are generated by the greedy algorithm producing the reduced GPT-PINNs of sizes $[2, 1, 1]$ to $[2, 9, 1]$. The GPT-PINNs are trained at the same set of collocation points as the SA-PINN  but with a learning rate of $0.0025$ and $2000$ epochs. 

\begin{figure}[!htbp]
    \centering
    \includegraphics[width=0.32\textwidth]{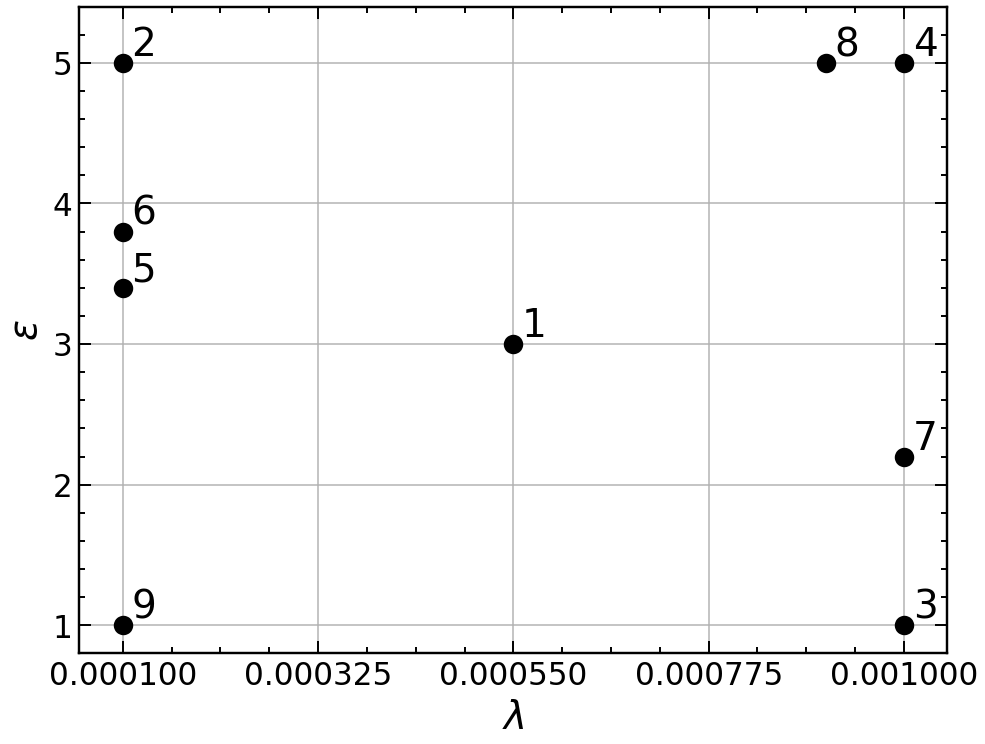}
    \includegraphics[width=0.32\textwidth]{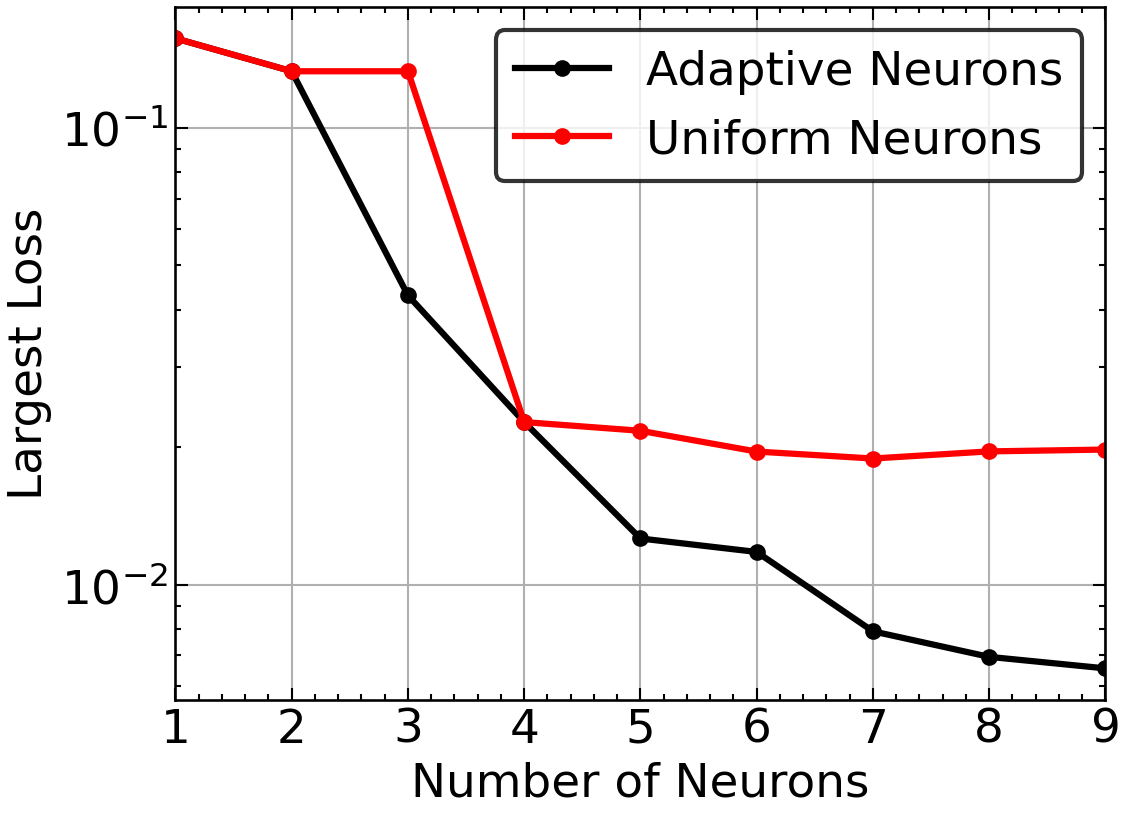}
    \includegraphics[width=0.32\textwidth]{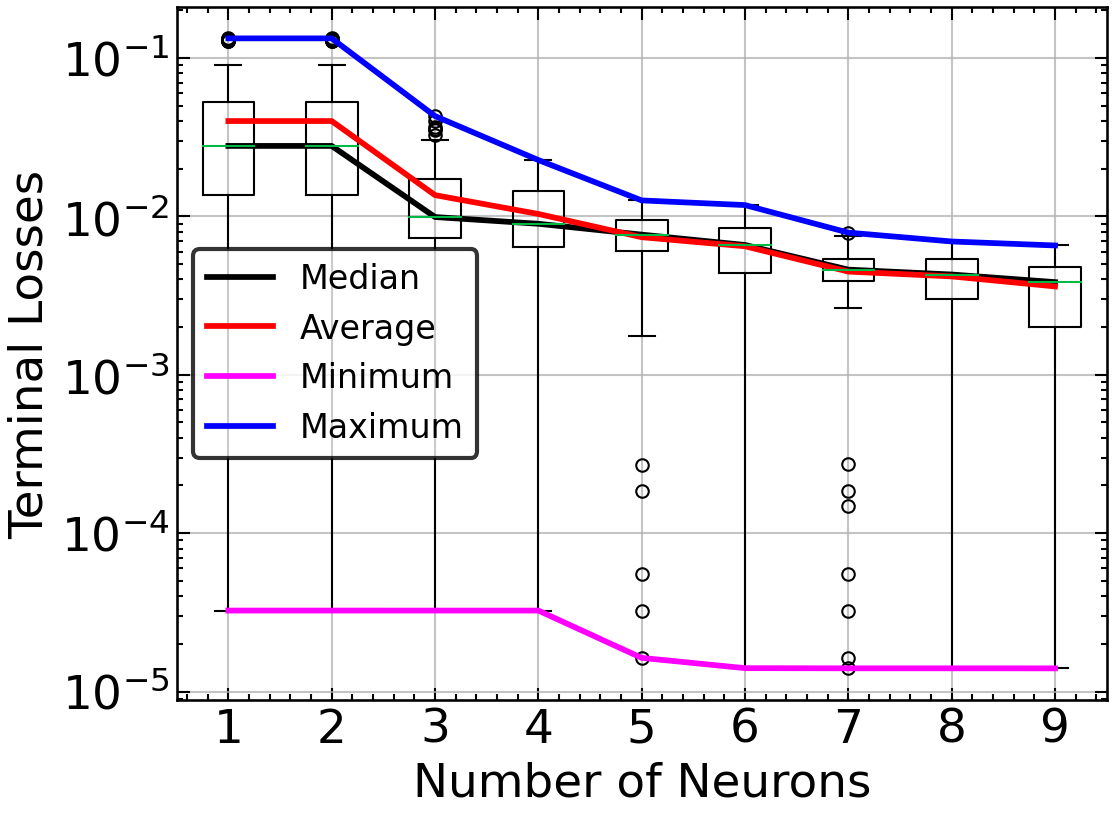}
    \caption{Allen-Cahn Equation training: The chosen parameter values (Left), worst-case GPT-PINN training
    losses (Middle), and the Box and Whisker plot of all GPT-PINN training losses (Right) during the outer-layer greedy
    training}
    \label{fig:pinnloss_AC}
\end{figure}
\begin{figure}[!htbp]
    \centering
    \includegraphics[width=0.32\textwidth]{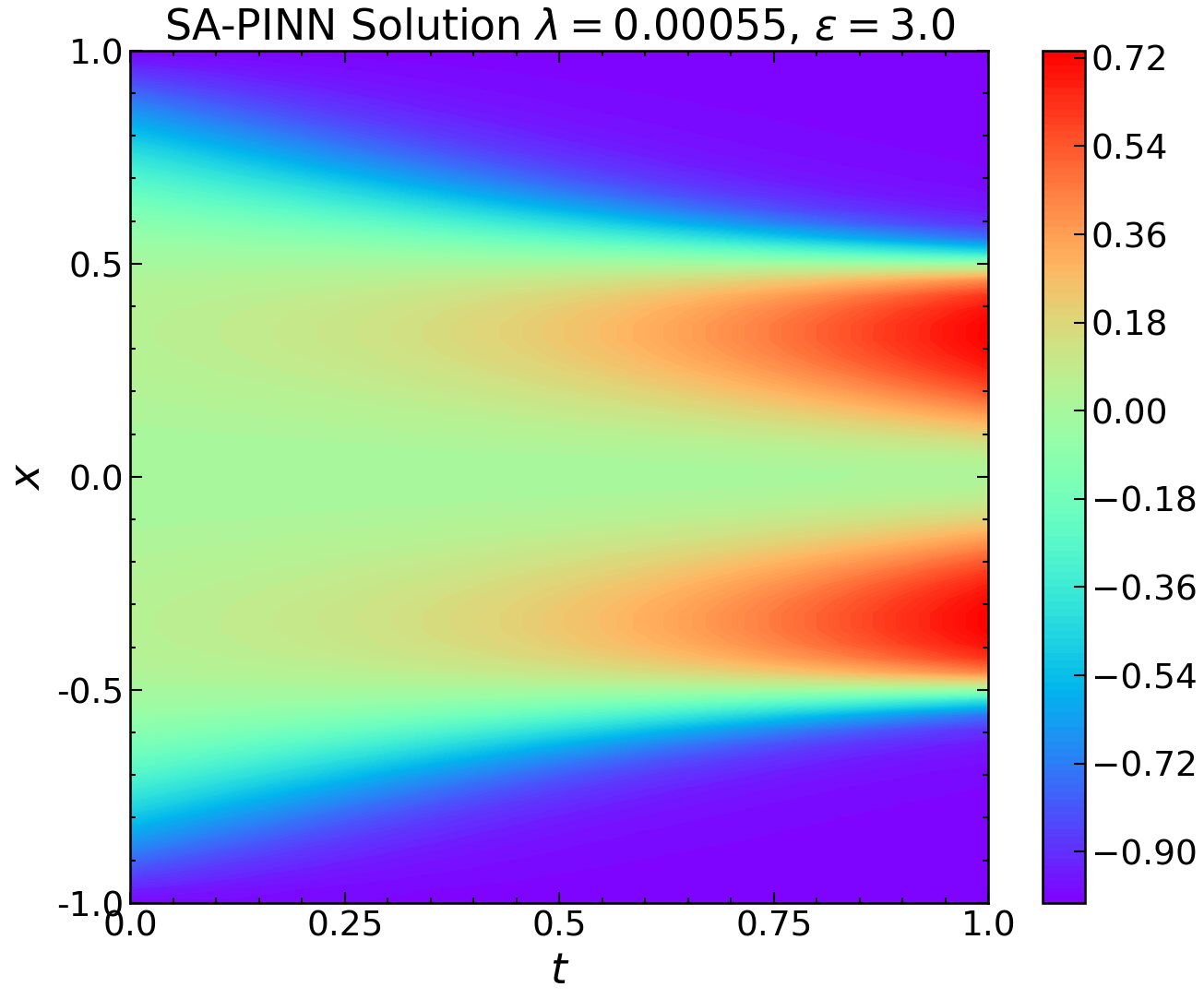}
    \includegraphics[width=0.32\textwidth]{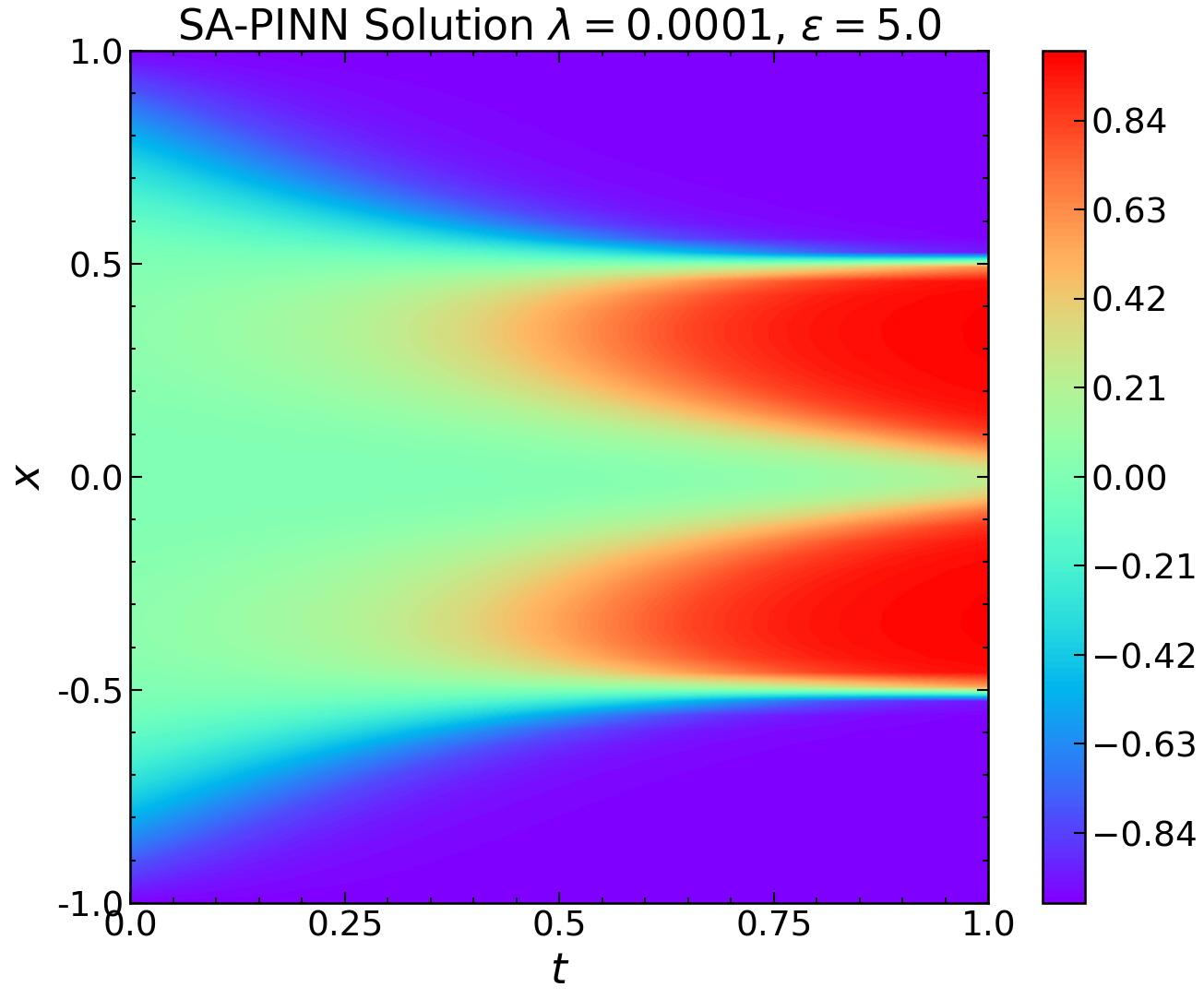}
    \includegraphics[width=0.32\textwidth]{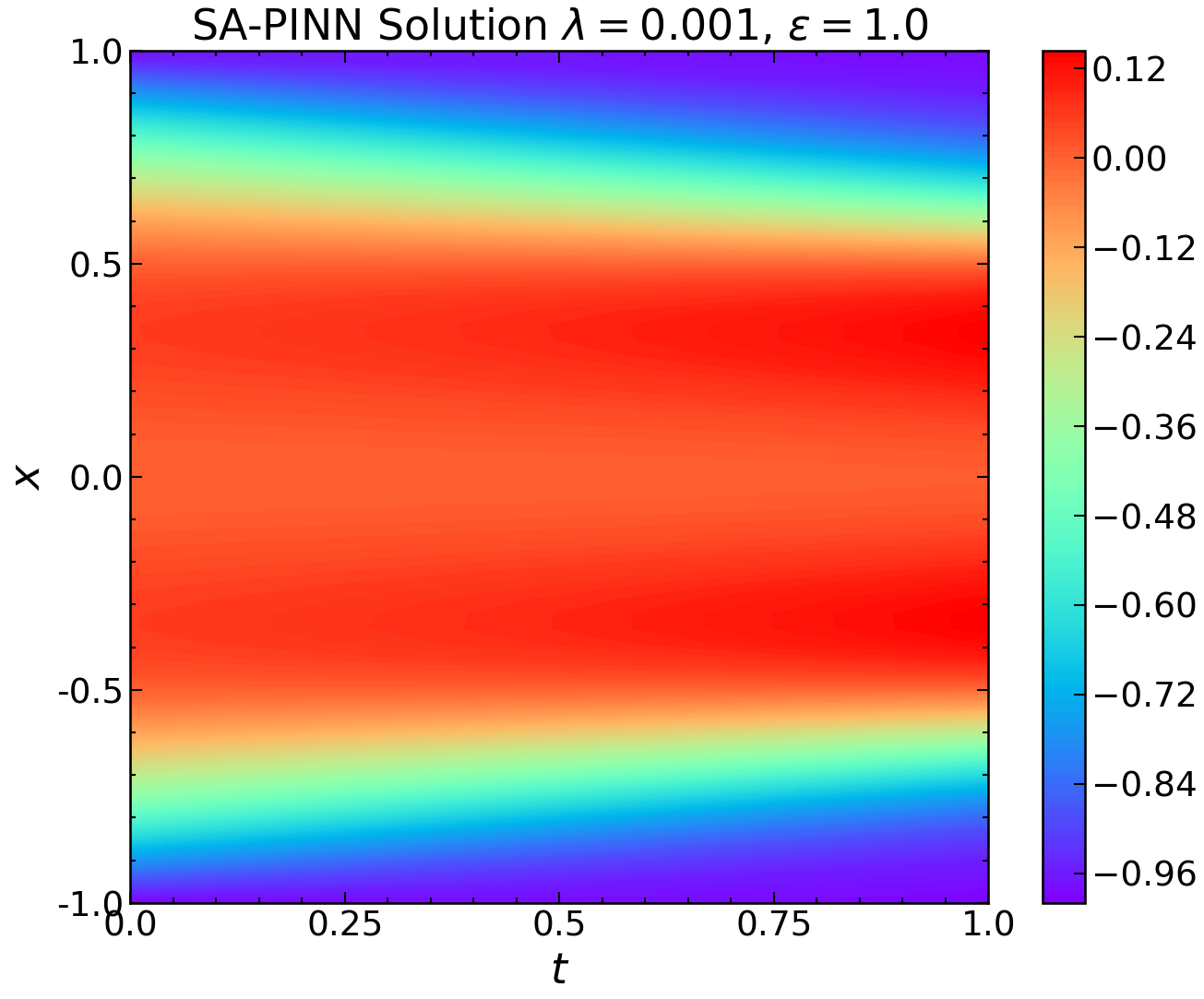}
    \caption{Allen-Cahn Equation: First three SA-PINN solutions found by the GPT-PINN that are used as the activation functions.}
    \label{fig:fullpinn_ac_sol}
\end{figure}

The GPT-PINN generates 9 neurons, i.e. SA-PINNs at $\{(\epsilon_i, \lambda_i)\}_{i=1}^{9}$. These parameter values and the worse-case offline training loss $\mathcal{L}_{\text{PINN}}^{\text{GPT}}(\bc(\bmu))$ after 2000 epochs as we increase the number of neurons (i.e. size of $\bc(\bmu)$) in the hidden layer of GPT-PINN are shown in Figure \ref{fig:pinnloss_AC}. Figure \ref{fig:fullpinn_ac_sol} shows the first three PINN solutions adaptively selected by GPT-PINN. 

\begin{figure}[!htbp]
    \centering
    \includegraphics[width=0.32\textwidth]{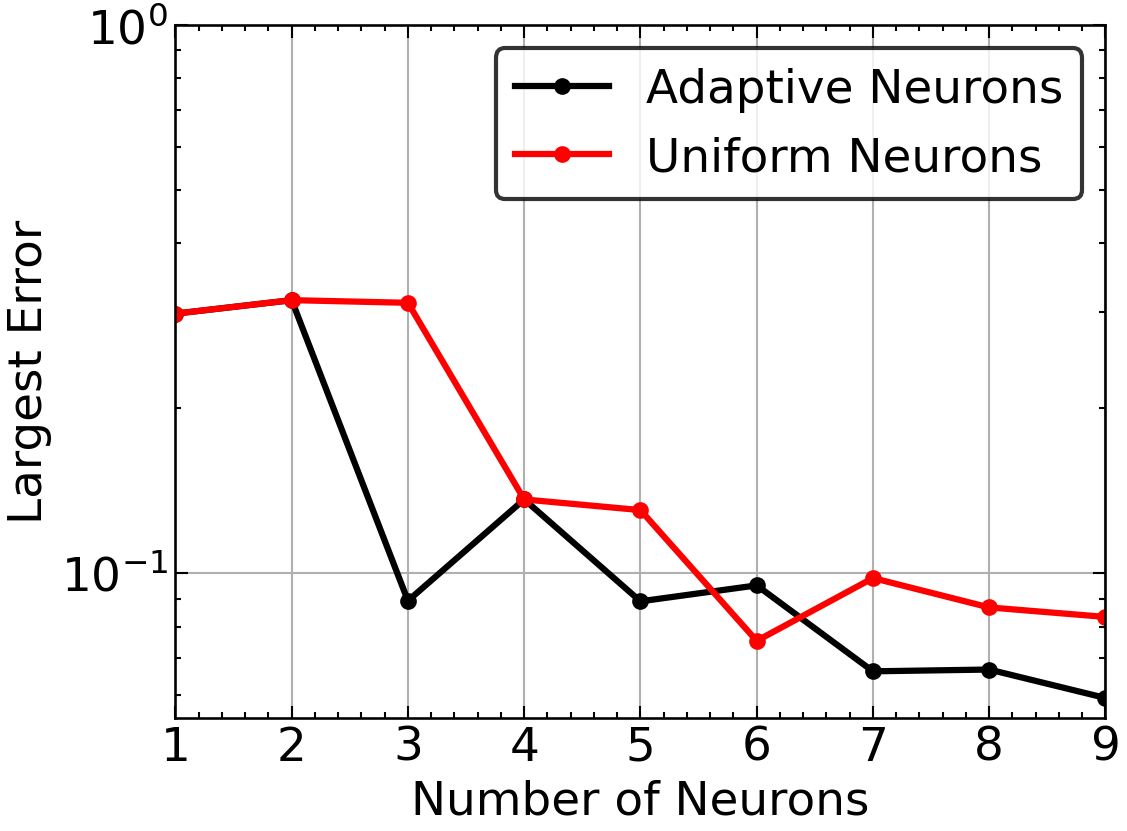}
    \includegraphics[width=0.32\textwidth]{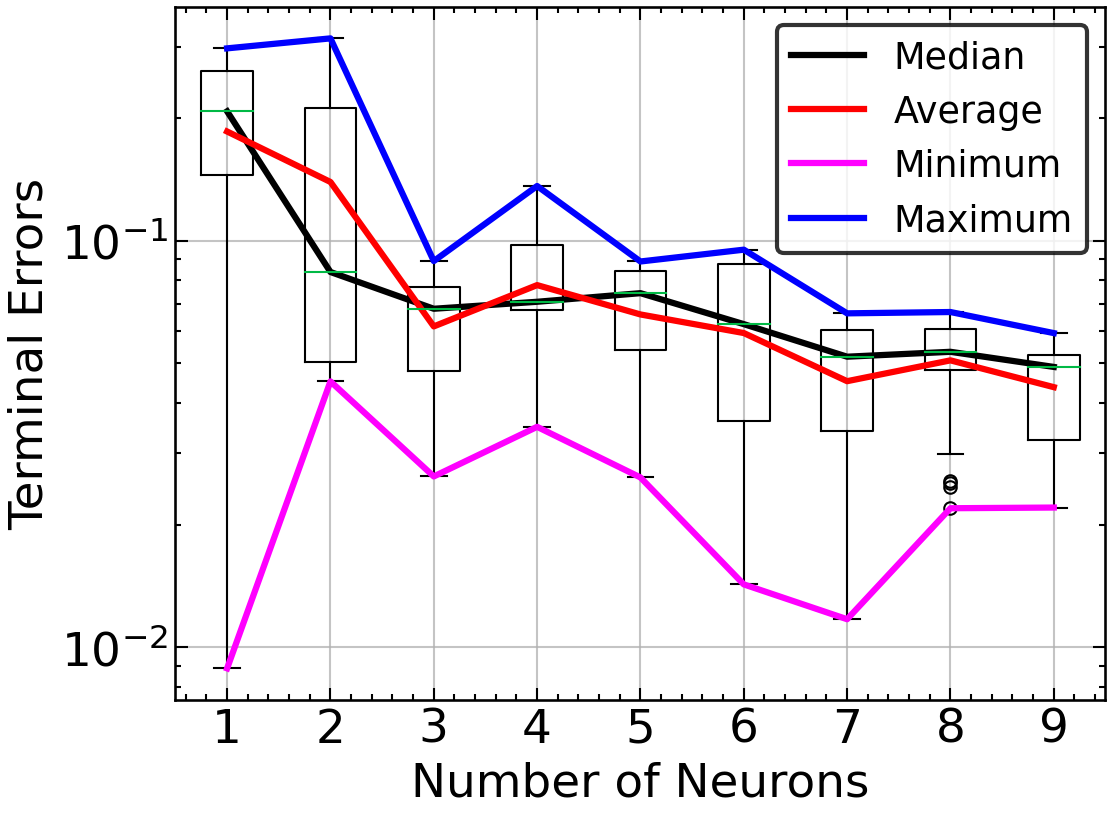}
    \includegraphics[width=0.32\textwidth]{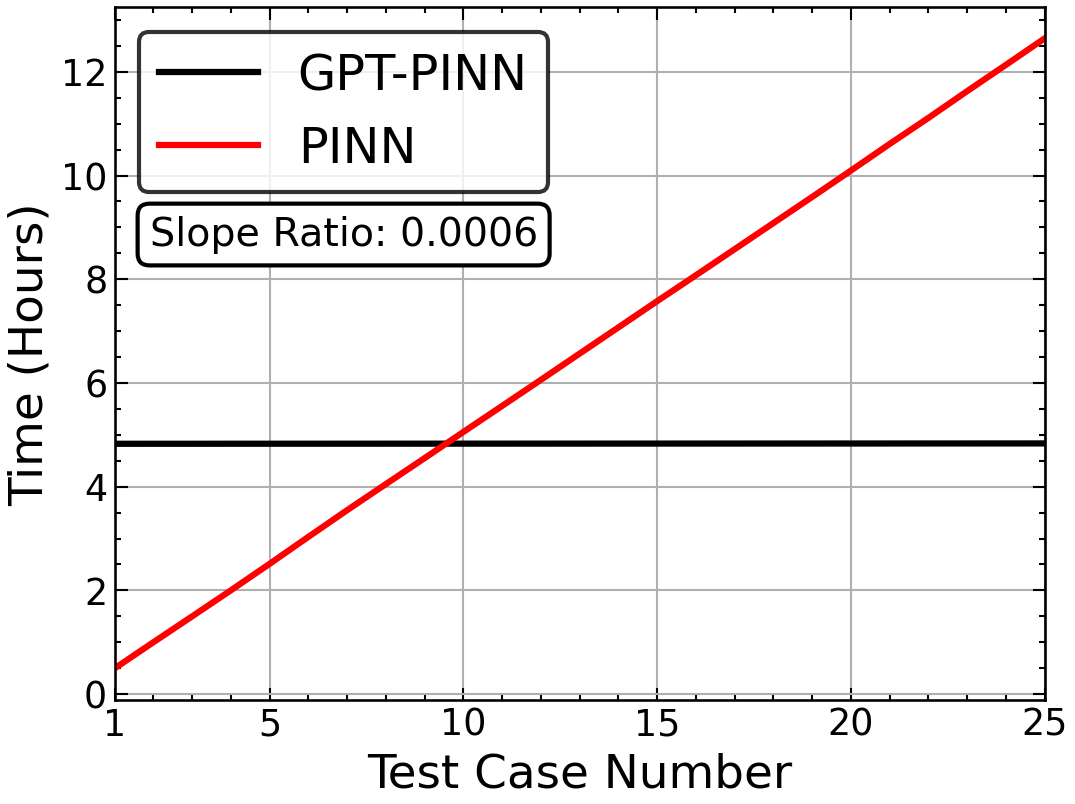}
    \caption{Allen-Cahn Equation testing: Worst-case test error of the GPT-PINN of various sizes (Left), Box and Whisker plot of all (Middle), and cumulative run time of the full PINN versus the GPT-PINN (Right)}
    \label{fig:pinnerror_AC}
\end{figure}
\begin{figure}[!htbp]
    \centering
    \includegraphics[width=0.32\textwidth]{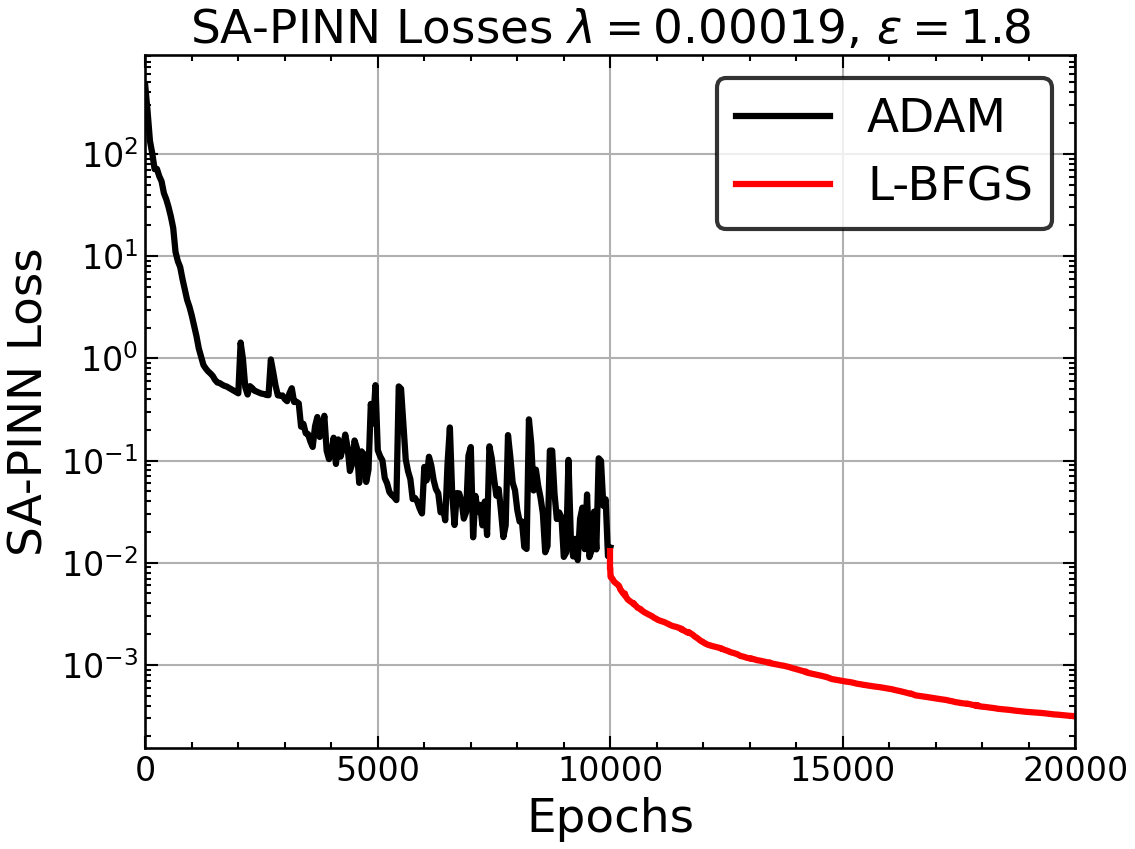}
    \includegraphics[width=0.32\textwidth]{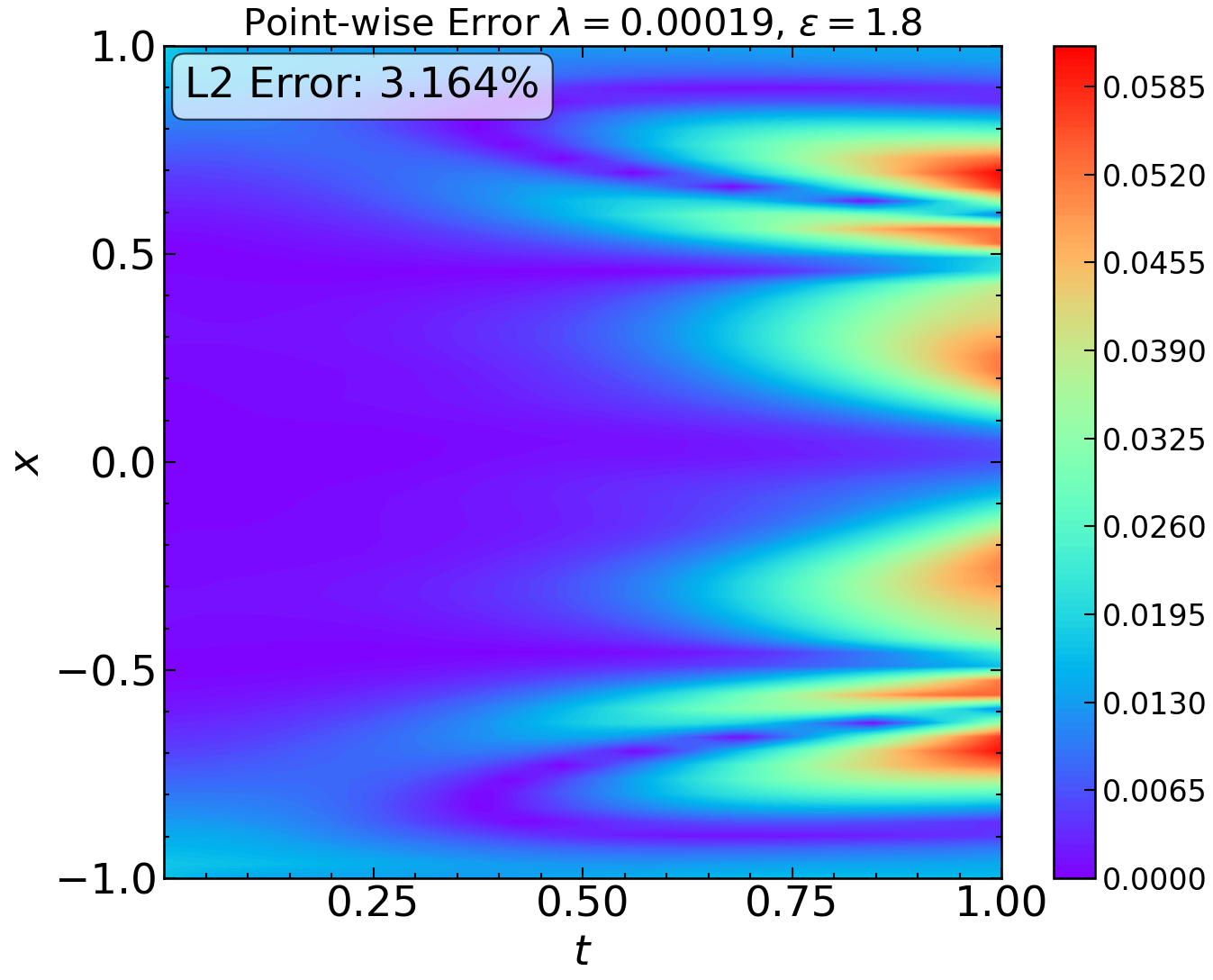}
    \includegraphics[width=0.32\textwidth]{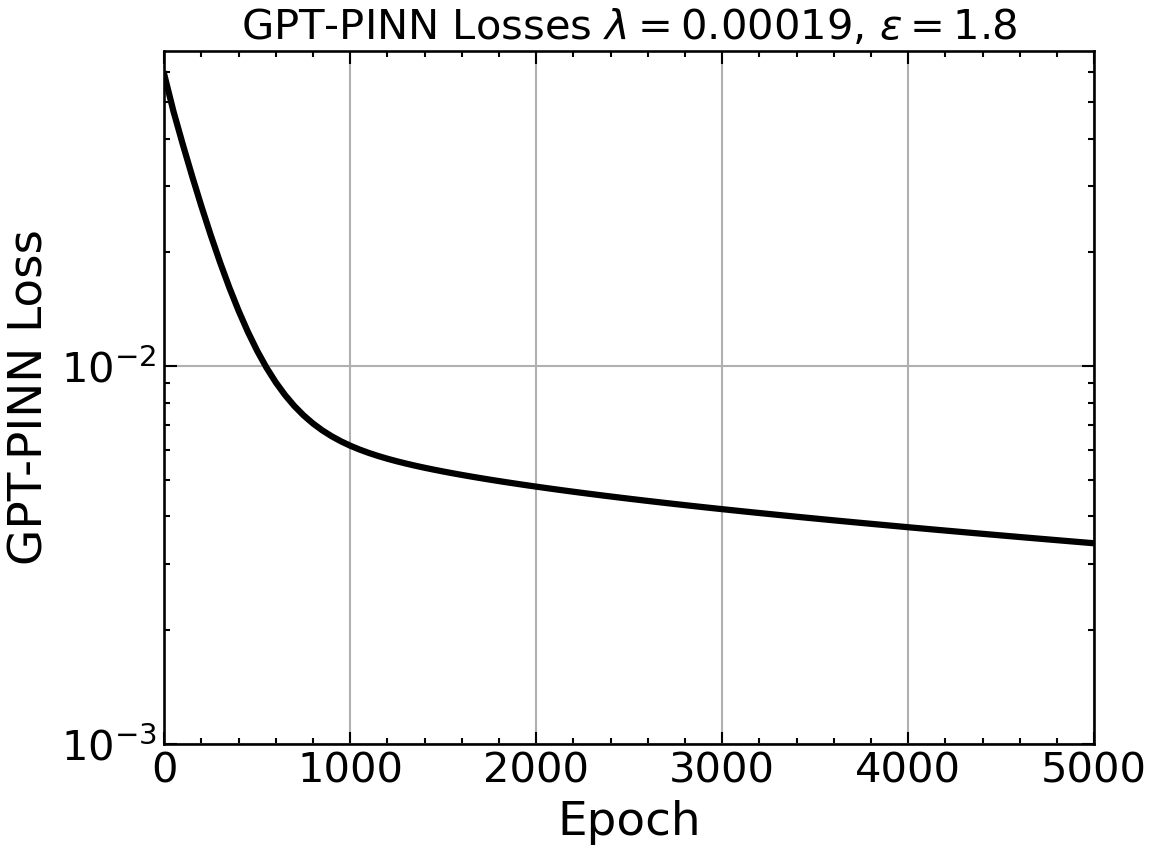}
    \includegraphics[width=0.32\textwidth]{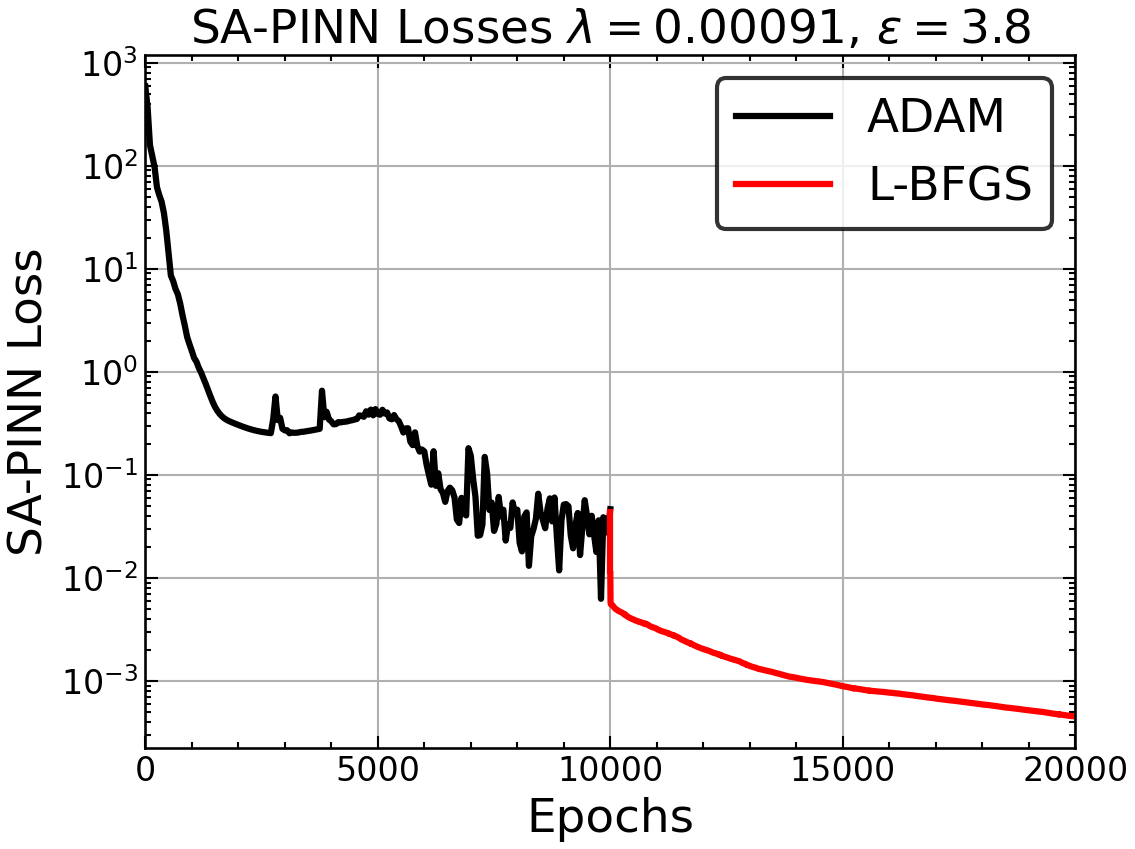}
    \includegraphics[width=0.32\textwidth]{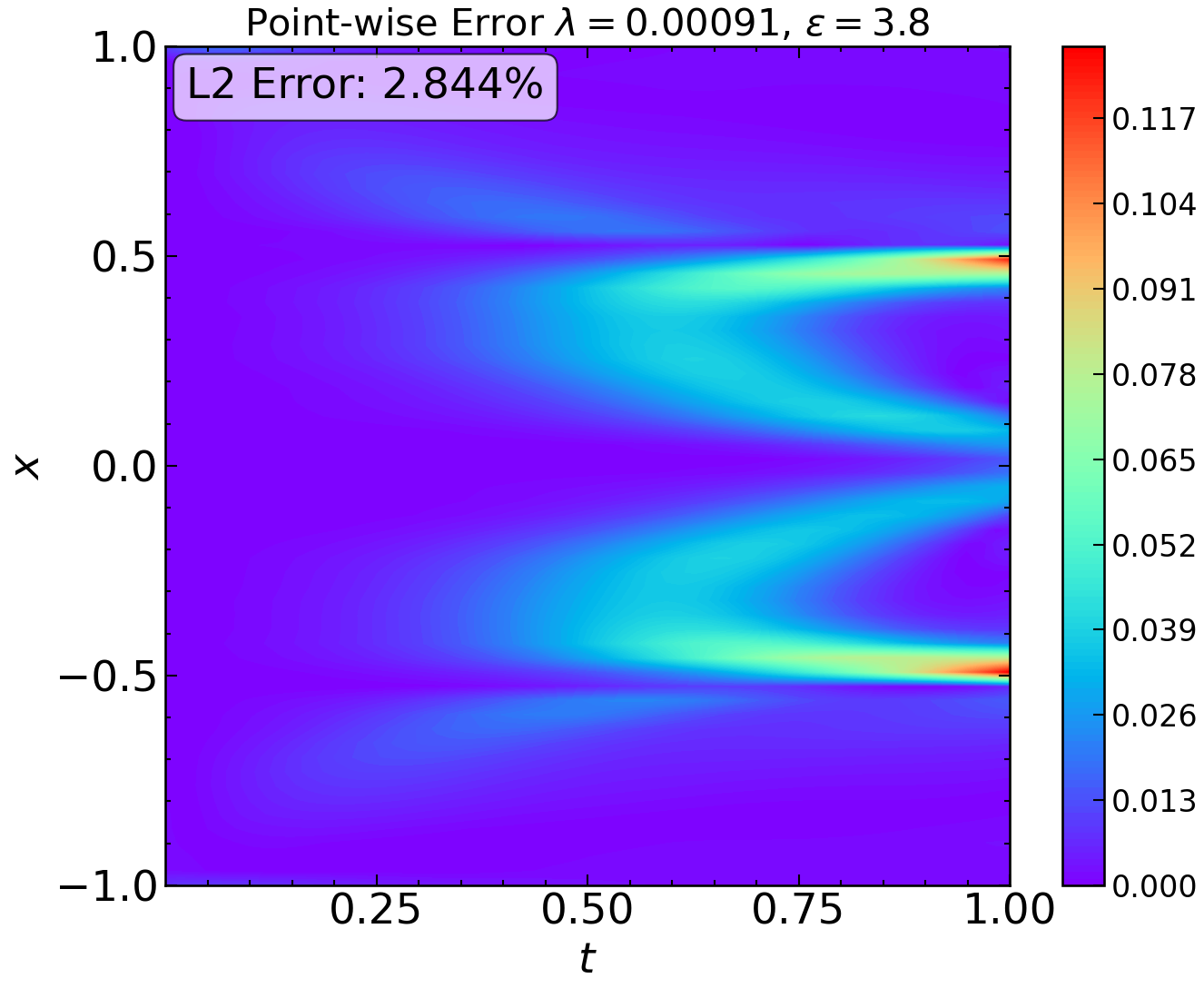}
    \includegraphics[width=0.32\textwidth]{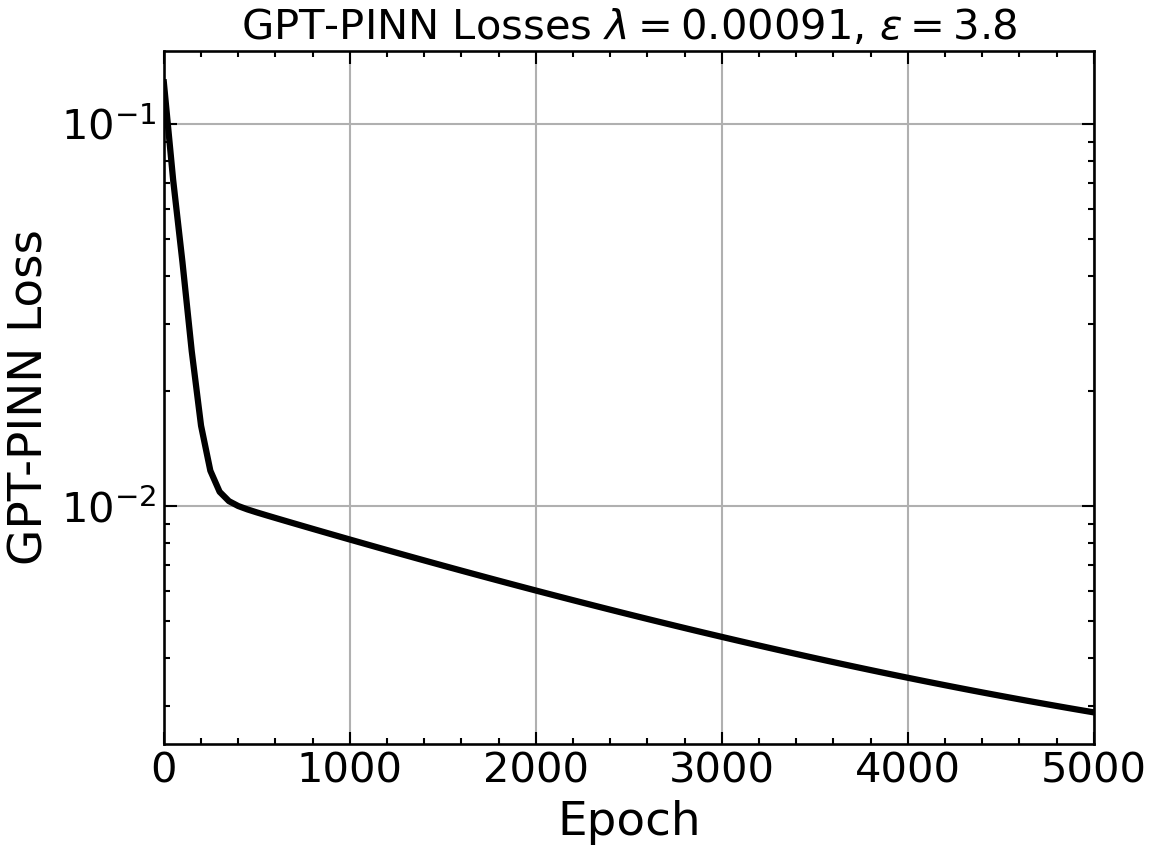}
    \caption{Allen-Cahn Equation: SA-PINN training loss (Left) and GPT-PINN training loss (Right) as functions of
    epochs for various parameters. Plotted in the middle are the point-wise errors of the corresponding GPT-PINN
    solution.}
    \label{fig:loss_v_epoch_ac}
\end{figure}

Next, we test the GPT-PINN on 25 parameter values. Figure \ref{fig:pinnerror_AC} displays the largest error for each size of the GPT-PINN and the cumulative run-time when both the SA-PINN and the GPT-PINN are repeatedly called. It is clear that the GPT-PINN line increases very slowly (at a relative speed of $0.0006$) and that it is again worthwhile to invest in GPT-PINN for a very modest number (9-10) of queries. We show the training losses as functions of epochs in Figure \ref{fig:loss_v_epoch_ac} for both the SA-PINN and GPT-PINN, with the latter again decaying more smoothly.

\section{Conclusion}
\label{sec:conclusion}

The proposed Generative Pre-Trained PINN (GPT-PINN)  is shown to mitigate two challenges faced by PINNs in the setting of parametric PDEs, namely the cost of training and over-parameterization. Being a hyper-reduced network with activation functions pre-trained full PINNs, GPT-PINN represents a brand-new meta-learning paradigm for parametric systems. 
With two main novelties, the design of network architecture including its special activation functions and the adoption of the training loss of the meta-network as an error indicator, and via tests on three differential families of parametric equations, 
we have shown that encompassing a very small number of well-chosen networks can generate surrogate PINNs across the entire parameter domain accurately and efficiently.

\appendix
\section{Detailed gradient of loss function for the Klein-Gordon case GPT-PINN}
\label{sec:appendix}

With the GPT-PINN formulation and considering the types of boundary and initial conditions for the equation given by \cref{eq:kg}, the loss function \cref{eq:loss-online} becomes 
\begin{align*}
\begin{split}
\mathcal{L}&_{\text{PINN}}^{\text{GPT}}(\bc(\bmu)) = \frac{1}{|\cC_o^r|} \sum_{(\bx,t) \in \cC_o^r}\left\lVert \frac{\partial^2}{\partial t^2}\left (\sum_{i=1}^n c_i(\bmu)\Psi^{\theta^i}_{\mathsf{NN}}\right)(\bx, t) + \alpha\frac{\partial^2}{\partial x^2}\left (\sum_{i=1}^n c_i(\bmu)\Psi^{\theta^i}_{\mathsf{NN}}\right)(\bx, t)  + \right. \\
& \left. \beta \left (\sum_{i=1}^n c_i(\bmu)\Psi^{\theta^i}_{\mathsf{NN}}\right)(\bx, t) + \gamma \left (\sum_{i=1}^n c_i(\bmu)\Psi^{\theta^i}_{\mathsf{NN}}\right)^2(\bx, t) + \bx\cos{(t)} - \bx^2\cos^2{(t)} \right\rVert_2^2 \\
& + \frac{1}{|\cC_\partial^r|}  \sum_{{(\bx,t) \in {\cC_\partial^r}} } \left\lVert \sum_{i=1}^n c_i(\bmu) \Psi^{\theta^i}_{\mathsf{NN}}(\bx, t) - u(\bx,t)\right\rVert_2^2 \\
& + \frac{1}{|\cC_i^r|} \sum_{\bx \in \cC_i^r} \left\lVert \sum_{i=1}^n c_i(\bmu) \Psi^{\theta^i}_{\mathsf{NN}}(\bx, 0) - u(\bx,0) \right\rVert_2^2 +  \frac{1}{|\cC_i^r|} \sum_{\bx \in \cC_i^r} \left\lVert \frac{\partial}{\partial t}\left (\sum_{i=1}^n c_i(\bmu)\Psi^{\theta^i}_{\mathsf{NN}}\right)(\bx, 0) - u_t(\bx, 0)\right\rVert_2^2
\end{split}
\end{align*}
with given $u(\bx,t)$ for $(\bx, t)\in \cC_\partial^r$ and $u(\bx,0)$ and $u_t(\bx,0)$ when $\bx \in \cC_i^r$. The $m^{\rm th}$ component of $\nabla_\bc \mathcal{L}_{\text{PINN}}^{\text{GPT}}(\bc)$ needed for the GPT-PINN training \cref{eq:c-update} then reads:
\begin{equation*}
\begin{split}
    \frac{\partial \mathcal{L}_{\text{PINN}}^{\text{GPT}}(\bc)}{\partial c_m} & = \frac{2}{|\cC_o^r|}\sum_{(\bx,t) \in \cC_o^r}  \Bigg(\bigg(\sum_{i=1}^n\big(c_i P^i_{tt} + \alpha c_i P^i_{xx} + \beta c_i P^i\big) + \gamma\big(\sum_{i=1}^{n}c_i P^i \big)^2 + \bx\cos{(t)} - \bx^2\cos^2{(t)}\bigg) \\
    & \cdot \bigg(P^m_{tt} + \alpha P^m_{xx} + \beta  P^m + 2\gamma\big(\sum_{i=1}^{n}c_i P^i \big)P^m \bigg)\Bigg) + \frac{2}{|\cC_\partial^r|}\sum_{(\bx,t) \in \cC_\partial^r}\Bigg(\bigg(\sum_{i=1}^{n}c_iP^i - u(\bx,t)\bigg)P^m \Bigg)\\
    & + \frac{2}{|\cC_i^r|}\sum_{\bx \in \cC_i^r}\Bigg(\bigg( \sum_{i=1}^n c_i P^i - u(\bx,0)\bigg)P^m\Bigg) +  \frac{2}{|\cC_i^r|}\sum_{\bx \in \cC_i^r}\Bigg(\bigg( \sum_{i=1}^n c_i P^i_t - u_t(\bx,0)\bigg)P_t^m\Bigg)
\end{split}
\end{equation*}
for $m = 1, \dots, n$. Here, for shortness of notation, we denote $\Psi^{\theta^i}_{\mathsf{NN}}(\bx, t)$ by $P^i(\bx,t)$ and omit $(\bx,t)$. 
For every full PINN $P^i$ identified by GPT-PINN, we would then just need to store the values of 
\[
P^i(\cC_o^r \cup \cC_\partial^r \cup (\cC_i^r\times \{0\})), \quad P^i_{xx}(\cC_o^r), \quad  P^i_{tt}(\cC_o^r), \quad P_t^i(\cC_i^r\times \{0\})
\]
for efficient online GPT-PINN training step of \cref{eq:c-update}.

\end{document}